\documentclass[11pt]{amsart}
\usepackage[table]{xcolor}
\usepackage{a4wide, amsmath, amsthm, mathabx, amscd, amsfonts, amssymb, calligra,
	caption, enumerate, epsfig, float, graphics, graphicx, hyperref, IEEEtrantools,
	mathrsfs, subcaption, textcomp, tikz, tikz-cd, verbatim, xypic, ytableau}
\usepackage[T1]{fontenc}

\newtheorem{theorem}{Theorem}[section]
\newtheorem{proposition}[theorem]{Proposition}
\newtheorem{conjecture}[theorem]{Conjecture}

\newtheorem{lemma}[theorem]{Lemma}

\theoremstyle{definition}

\newtheorem{remark}[theorem]{Remark}
\newtheorem{example}[theorem]{Example}
\newtheorem{definition}[theorem]{Definition}

\def\Cl{\calligra{Cl}}

\def\cC{\mathcal{C}}

\def\hol{\operatorname{hol}}
\def\Gr{\operatorname{Gr}}
\def\vGr{\operatorname{Gr}^\vee}
\def\QH{\operatorname{QH}}
\def\Fuk{\mathcal{F}}

\def\hom{\operatorname{hom}}

\def\HF{\operatorname{HF}}
\newcommand{\vgeq}{\mathbin{\rotatebox[origin=c]{-90}{$\geq$}}}
\newcommand{\brgeq}{\mathbin{\rotatebox[origin=c]{-45}{$\geq$}}}
\newcommand{\trgeq}{\mathbin{\rotatebox[origin=c]{45}{$\geq$}}}

\def\D{{\bf D}}
\def\Cl{{\calligra Cl}}

\definecolor{darkgreen}{RGB}{0,153,0}
\definecolor{darkred}{RGB}{204,0,0}
\definecolor{darkblue}{RGB}{0,51,204}
\definecolor{red}{RGB}{242,43,29}
\hypersetup{colorlinks,linkcolor={darkblue},citecolor={darkblue},urlcolor={darkblue}}  

\begin{document}

\title{Fukaya category of Grassmannians: rectangles}
\author{Marco Castronovo}
\date{}
\address{Rutgers University - Hill Center for the Mathematical Sciences}
\email{marco.castronovo@rutgers.edu}

\thanks{Partially supported by NSF grant DMS 1711070. Any opinions, findings,
and conclusions or recommendations expressed in this material are those of the author and do not necessarily
reflect the views of the National Science Foundation.}

\begin{abstract}
We show that the monotone Lagrangian torus fiber of the Gelfand-Cetlin integrable system
on the complex Grassmannian $\operatorname{Gr}(k,n)$ supports generators
for all maximum modulus summands in the spectral decomposition of the Fukaya category over
$\mathbb{C}$, generalizing the example of the Clifford torus in projective space. We introduce
an action of the dihedral group $D_n$ on the Landau-Ginzburg mirror proposed by Marsh-Rietsch \cite{MR}
that makes it equivariant and use it to show that, given a lower modulus, the torus
supports nonzero objects in none or many summands of the Fukaya category with that modulus.
The alternative is controlled by the vanishing of rectangular Schur polynomials at the $n$-th
roots of unity, and for $n=p$ prime this suffices to give a complete set of generators and prove homological mirror
symmetry for $\Gr(k,p)$.
\end{abstract}

\maketitle
\thispagestyle{empty}
\tableofcontents

\section*{Introduction}\label{Intro}

According to a conjecture described from a mathematical viewpoint by Auroux \cite{AuT}, there should be a construction that starts
from the choice of an anti-canonical divisor $D\subset X$ in a compact K\"{a}hler manifold and a holomorphic volume
form $\Omega$ on $X\setminus D$, and produces a complex manifold
$X^\vee$, sometimes referred to as the \textit{Landau-Ginzburg mirror}, with a holomorphic function
$W:X^\vee\to \mathbb{C}$ called the \textit{superpotential}. The terminology is inspired by the work of string
theorists on dualities between \textit{D-branes}, see Hori-Iqbal-Vafa \cite{HIV}. Roughly, $X^\vee$ should arise as moduli space of Lagrangian
tori $L\subset X\setminus D$ equipped with rank one $\mathbb{C}$-linear local systems $\xi$ and calibrated
by $\operatorname{Re}\Omega$, while $W$ should be the obstruction to the Floer operator
squaring to zero in the space of Floer cochains
of $L$, essentially determined by the pseudo-holomorphic disks bounding $L$ in $X$. These disks carry
information about the symplectic topology of $L\subset X$: for instance, Vianna \cite{Vi} used them to prove that there
are infinitely many monotone Lagrangian tori in $\mathbb{P}^2$ not Hamiltonian isotopic to each other.
\par
One way to make precise the idea of studying all Lagrangians $L\subset X$ together has
been described in Fukaya-Oh-Ohta-Ono \cite{FOOO}, and
consists in constructing some flavor of an $A_\infty$-category $\mathcal{F}(X)$ whose objects
are Lagrangians, morphisms are Floer cochains, and structure maps count pseudo-holomorphic disks.
We shall use a variant of Seidel's construction \cite{Se} of this category in the exact setting,
which is described by Sheridan \cite{Sh} and works for compact \textit{monotone} manifolds. This
case is relevant for us because it includes Fano smooth projective varieties over $\mathbb{C}$. We
briefly summarize this framework in the Setup section. The key structural property that we will
use is the \textit{spectral decomposition} of the Fukaya category
$$\Fuk(X)=\bigoplus_{\lambda}\Fuk_\lambda(X)$$
in summands labelled by the eigenvalues of the operator $c_1\star$ of multiplication by the first
Chern class acting on quantum cohomology $\QH(X)$.
\par
As suggested by Kontsevich \cite{Ko} one can consider the \textit{derived Fukaya category}
$\D\Fuk(X)$, a triangulated category, and phrase constructions like
the one mentioned at the beginning in terms of \textit{equivalences} with triangulated
categories carrying informations about sheaves on $X^\vee$. In our case, the relevant partner
for $\D\Fuk_\lambda(X)$ will be the \textit{derived category of singularities} $\D\mathcal{S}(W^{-1}(\lambda))$
introduced by Orlov \cite{Or}, which measures to what extent coherent sheaves on the fiber
$W^{-1}(\lambda)$ fail to have finite resolutions by locally free sheaves. With these tools in place, \textit{homological mirror
symmetry} holds if one can establish equivalences of triangulated categories
$$\D\Fuk_\lambda(X)\simeq \D\mathcal{S}(W^{-1}(\lambda))$$
for all eigenvalues $\lambda$.
\par
A question one could start with is to find sets of \textit{generators} for these triangulated
categories. On one side Dyckerhoff \cite{Dy} showed that, whenever $W$ has isolated
singularities, the category of singularities is generated by skyscraper sheaves
at the singular points. On the symplectic side no such general statement is known, and generators
have been described only in special cases. For $X=\mathbb{P}^n$ Cho \cite{Ch} showed
that the Clifford Lagrangian torus supports $n+1$ local
systems with nonzero Floer cohomology, corresponding to the critical points of $W:X^\vee\to \mathbb{C}$
with $X^\vee=(\mathbb{C}^*)^n$ and
$$W = z_1 + \ldots + z_n + \frac{1}{z_1\cdots z_n}\quad .$$
This picture generalizes to arbitrary
$X=X(\Delta)$ smooth projective Fano toric varieties over $\mathbb{C}$, where the Clifford torus
is replaced by the monotone Lagrangian torus fiber over the barycenter of the polytope $\Delta$
and again $X^\vee=(\mathbb{C}^*)^n$, with $W$ determined by $\Delta$
and now one has $\chi(X(\Delta))$ local systems corresponding to its critical points. See the work
of Fukaya-Oh-Ohta-Ono \cite{FOOOt} (also relevant for more general settings).
\par
\begin{figure}[H]
  \centering
   %add desired spacing between images, e. g. ~, \quad, \qquad, \hfill etc. 
    %(or a blank line to force the subfigure onto a new line)
    \begin{subfigure}[b]{0.3\textwidth}
        \includegraphics[width=\textwidth]{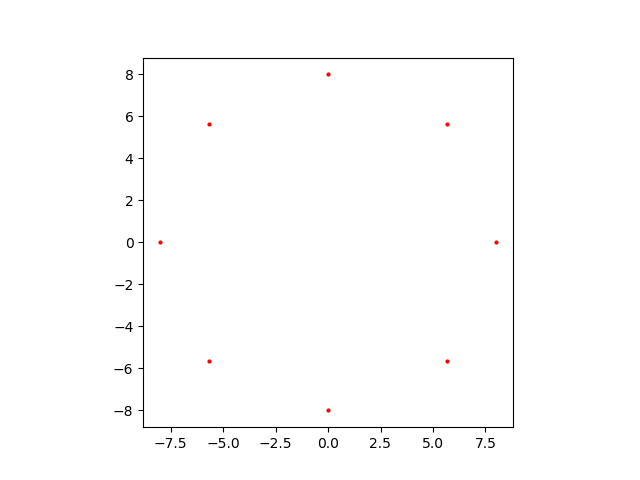}
        \caption{$\Gr(1,8)$}
        \label{fig:flower(1,8)}
    \end{subfigure}
    \begin{subfigure}[b]{0.3\textwidth}
        \includegraphics[width=\textwidth]{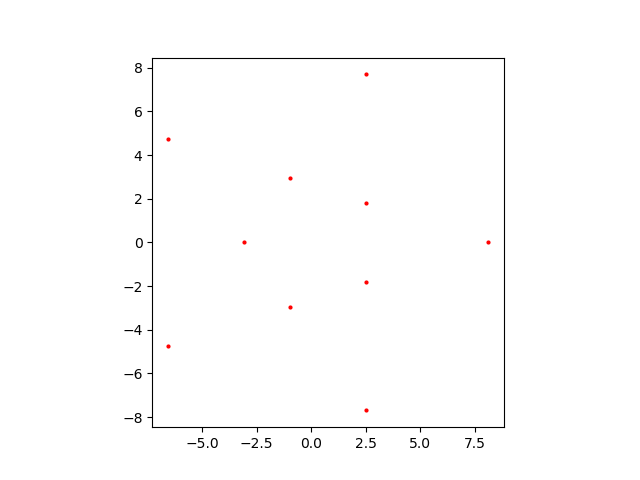}
        \caption{$\Gr(2,5)$}
        \label{fig:flower(2,5)}
    \end{subfigure}
     %add desired spacing between images, e. g. ~, \quad, \qquad, \hfill etc. 
    %(or a blank line to force the subfigure onto a new line)
    \begin{subfigure}[b]{0.3\textwidth}
        \includegraphics[width=\textwidth]{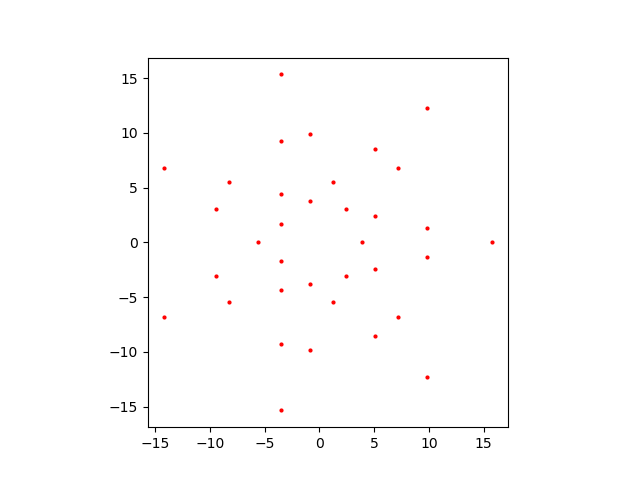}
        \caption{$\Gr(3,7)$}
        \label{fig:flower(3,7)}
    \end{subfigure}
    
     %add desired spacing between images, e. g. ~, \quad, \qquad, \hfill etc. 
      %(or a blank line to force the subfigure onto a new line)
    \begin{subfigure}[b]{0.3\textwidth}
        \includegraphics[width=\textwidth]{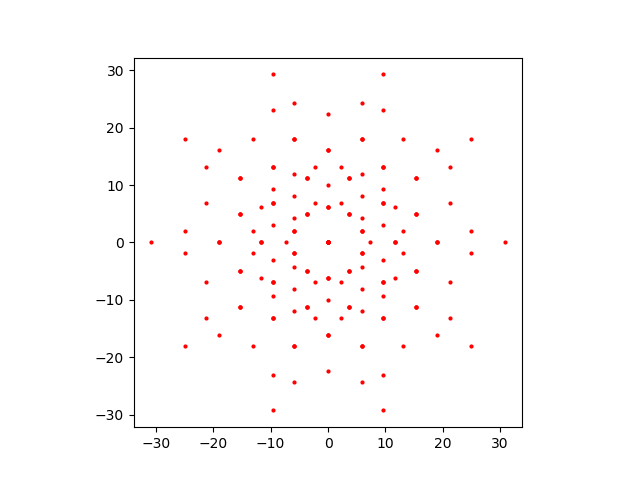}
        \caption{$\Gr(4,10)$}
        \label{fig:flower(4,10)}
    \end{subfigure}
    %add desired spacing between images, e. g. ~, \quad, \qquad, \hfill etc. 
    %(or a blank line to force the subfigure onto a new line)
    \begin{subfigure}[b]{0.3\textwidth}
        \includegraphics[width=\textwidth]{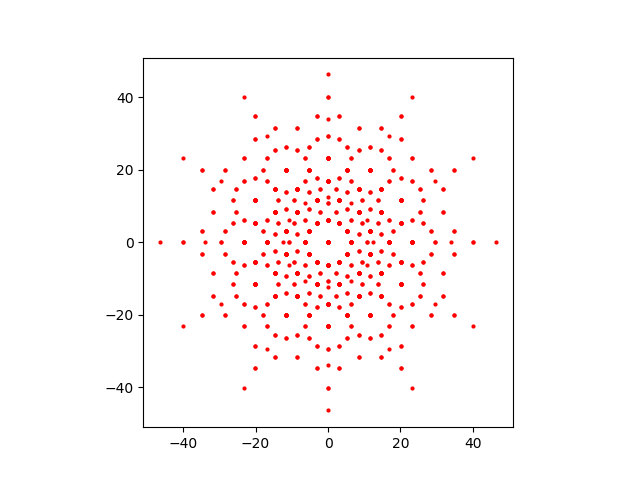}
        \caption{$\Gr(6,12)$}
        \label{fig:flower(6,12)}
    \end{subfigure}
    %add desired spacing between images, e. g. ~, \quad, \qquad, \hfill etc. 
    %(or a blank line to force the subfigure onto a new line)
    \begin{subfigure}[b]{0.3\textwidth}
        \includegraphics[width=\textwidth]{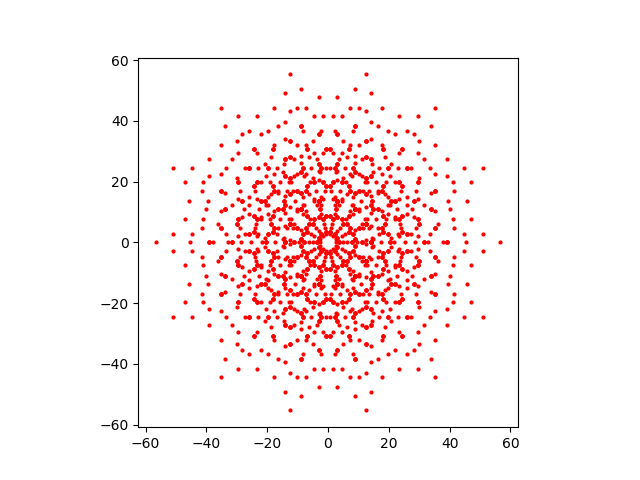}
        \caption{$\Gr(5,14)$}
        \label{fig:flower(5,14)}
    \end{subfigure}

    \caption{Eigenvalues of $c_1\star$ acting on $\QH(\Gr(k,n))$.}
    \label{fig:flowers}
\end{figure}

A natural case to consider next is the one of Grassmannians $X=\Gr(k,n)$ parametrizing
$k$-dimensional linear subspaces of $\mathbb{C}^n$. This exhibits some novel features that are not apparent
in the toric case: building on Peterson's general presentation of the quantum cohomology of flag
varieties, Marsh-Rietsch \cite{MR} proposed a Landau-Ginzburg mirror $X^\vee=\vGr(k,n)$ that is
not a single complex torus, but rather a glueing of complex torus charts. We review this in Section
\ref{Section1}, where we show (Proposition \ref{Prop1}) how
certain diagrams with \textit{dihedral symmetry} index the summands of the Fukaya category: see
Figure \ref{fig:flowers}.

Section \ref{Section2} precises what we mean by complex torus chart of the Landau-Ginzburg mirror,
focusing on one that we call \textit{rectangular chart}. We give a criterion (Proposition \ref{Prop2}) for deciding when a
critical point of $W$ belongs to the rectangular chart, formulated in terms of vanishing of certain Schur
polynomials at roots of unity.
\par
The Grassmannian $\Gr(k,n)$ supports an analogue of the toric moment map,
which is an integrable system introduced by Guillemin-Sternberg \cite{GS} and
called \textit{Gelfand-Cetlin} integrable system.
This is the topic of Section \ref{Section3}, where we use a \textit{toric degeneration} argument of
Nishinou-Nohara-Ueda \cite{NNU} and a combinatorial description of the faces of its
image polytope due to An-Cho-Kim \cite{ACK} to write down a formula for the Maslov 2 disk potential
of a monotone Lagrangian torus fiber of the system (Proposition \ref{Prop3}).
We think of this Lagrangian as mirror to the rectangular chart of $\vGr(k,n)$.
\par
The aim of this paper is to focus on local systems supported on this torus, that we will
call from now on Gelfand-Cetlin torus; it generalizes the monotone
Clifford torus in $\mathbb{P}^n$ to the other Grassmannians, and
we want to understand to what extent one can use objects supported on it to generate the Fukaya
category.
\par
In order to state the main theorems, we introduce now some notation. Call $T^{k(n-k)}\subset \Gr(k,n)$
the Gelfand-Cetlin torus, and $\gamma_{ij}\in H_1(T^{k(n-k)};\mathbb{Z})$ for $1\leq i\leq k,\, 1\leq j\leq n-k$
the basis of cycles induced by the integrable system. For each set $I$ of $k$
distinct roots of $x^n=(-1)^{k+1}$, define a local system on $T^{k(n-k)}$ whose holonomy
is given by
$$\operatorname{hol}_I(\gamma_{ij})=\frac{S_{(k+1-i)\times j}(I)}{S_{(k-i)\times (j-1)}(I)} \quad .$$
In this formula, $S_{p\times q}$ is the $k$-variables Schur polynomial of a rectangular Young diagram
$p\times q$ in a $k\times (n-k)$ grid. The definition above makes sense only when the denominator
is nonzero, in which case $T^{k(n-k)}_I$ denotes the corresponding object
of the Fukaya category. When the denominator is zero, the object $T^{k(n-k)}_I$ is not defined.
Also consider the dihedral group
$$D_n = \langle r,\, s\, |\, r^n=s^2=1,\, rs=sr^{-1} \rangle$$
and its action on the sets $I$ via $rI=e^{2\pi i/n}I$ and $sI=\overline{I}$.

\begin{theorem}\label{thm:main-thm-1-statement}
Choosing $I_0$ to be the set of $k$ roots of $x^n=(-1)^{k+1}$ closest to 1, the objects
obtained by giving $T^{k(n-k)}$ the different local systems
$$T^{k(n-k)}_{I_0},\; T^{k(n-k)}_{rI_0},\; T^{k(n-k)}_{r^2I_0},\; \ldots \;,\; T^{k(n-k)}_{r^{n-1}I_0}$$
are defined and split-generate the $n$ summands $\Fuk_\lambda(\Gr(k,n))$ of the monotone
Fukaya category with maximum $|\lambda|$.
\end{theorem}

The proof of this theorem is based on Sheridan's extension \cite{Sh} of the generation
criterion of Abouzaid \cite{Ab}. The main contribution is the observation that, under
an explicit open embedding of schemes
$$\theta_R:(\mathbb{C}^\times)^{k(n-k)}\to \vGr(k,n)$$
that identifies the space of local systems on the Gelfand-Cetlin torus with the rectangular
chart of the mirror, the superpotential $W$ pulls back to the
Maslov 2 disk potential of the Gelfand-Cetlin torus. One can generate each maximum modulus
summand of the Fukaya category with a single nonzero object, and
such objects are found by endowing the Gelfand-Cetlin torus with local systems whose
holonomies are rotations of the unique critical point of $W$ lying in the
\textit{totally positive} part $\vGr(k,n)_{\geq 0}\subset\vGr(k,n)$ of the mirror.
This is the locus where the classic Pl\"{u}cker coordinates are real and positive,
and its properties have been the focus of several works in representation theory,
combinatorics, topology and mirror symmetry: see \cite{Lu}, \cite{Po}, \cite{KLS},
\cite{STWZ}, \cite{RW}.
\par
\begin{figure}[H]
  \centering
   %add desired spacing between images, e. g. ~, \quad, \qquad, \hfill etc. 
    %(or a blank line to force the subfigure onto a new line)
    \begin{subfigure}[b]{0.3\textwidth}
        \includegraphics[width=\textwidth]{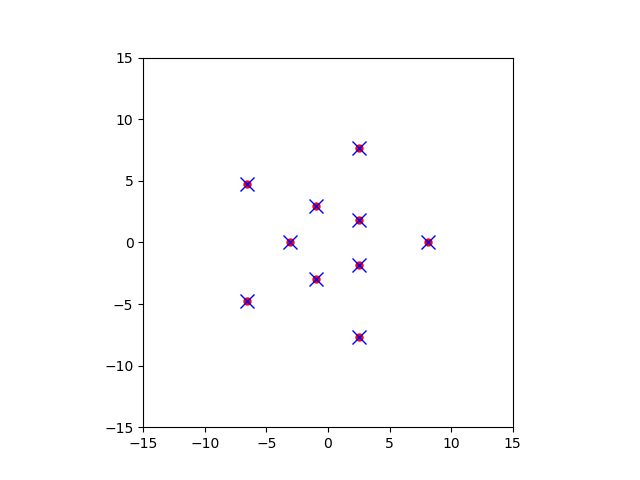}
        \caption{$\Gr(2,5)$}
        \label{fig:branes(2,5)}
    \end{subfigure}
    \begin{subfigure}[b]{0.3\textwidth}
        \includegraphics[width=\textwidth]{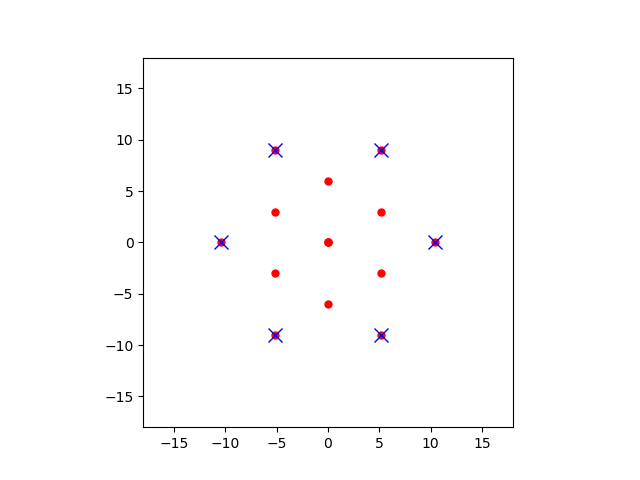}
        \caption{$\Gr(2,6)$}
        \label{fig:branes(2,6)}
    \end{subfigure}
     %add desired spacing between images, e. g. ~, \quad, \qquad, \hfill etc. 
    %(or a blank line to force the subfigure onto a new line)
    \begin{subfigure}[b]{0.3\textwidth}
        \includegraphics[width=\textwidth]{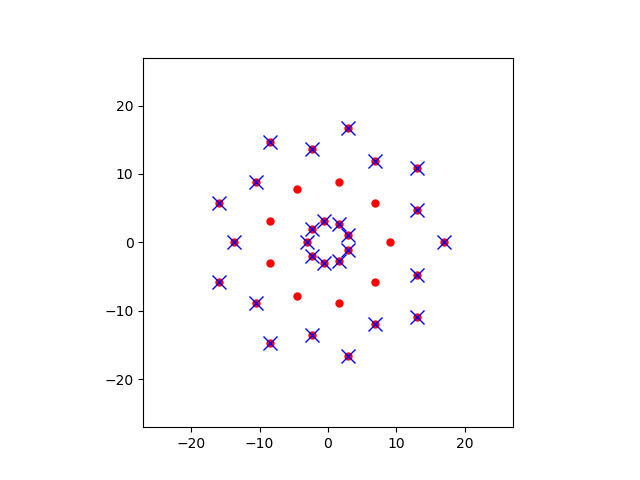}
        \caption{$\Gr(2,9)$}
        \label{fig:branes(2,9)}
    \end{subfigure}
    
     %add desired spacing between images, e. g. ~, \quad, \qquad, \hfill etc. 
      %(or a blank line to force the subfigure onto a new line)
    \begin{subfigure}[b]{0.3\textwidth}
        \includegraphics[width=\textwidth]{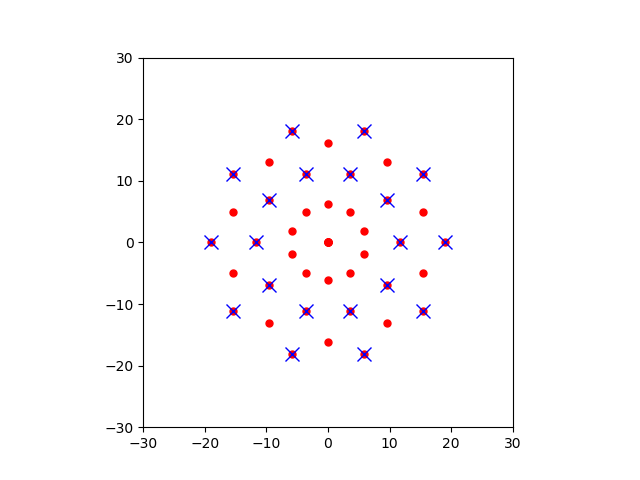}
        \caption{$\Gr(2,10)$}
        \label{fig:branes(2,10)}
    \end{subfigure}
    %add desired spacing between images, e. g. ~, \quad, \qquad, \hfill etc. 
    %(or a blank line to force the subfigure onto a new line)
    \begin{subfigure}[b]{0.3\textwidth}
        \includegraphics[width=\textwidth]{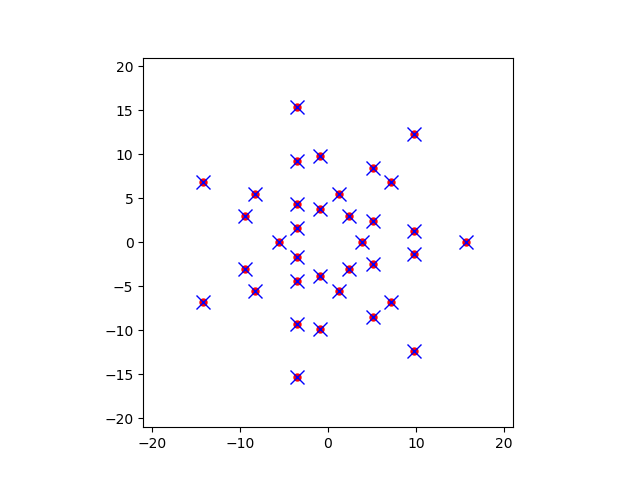}
        \caption{$\Gr(3,7)$}
        \label{fig:branes(3,7)}
    \end{subfigure}
    %add desired spacing between images, e. g. ~, \quad, \qquad, \hfill etc. 
    %(or a blank line to force the subfigure onto a new line)
    \begin{subfigure}[b]{0.3\textwidth}
        \includegraphics[width=\textwidth]{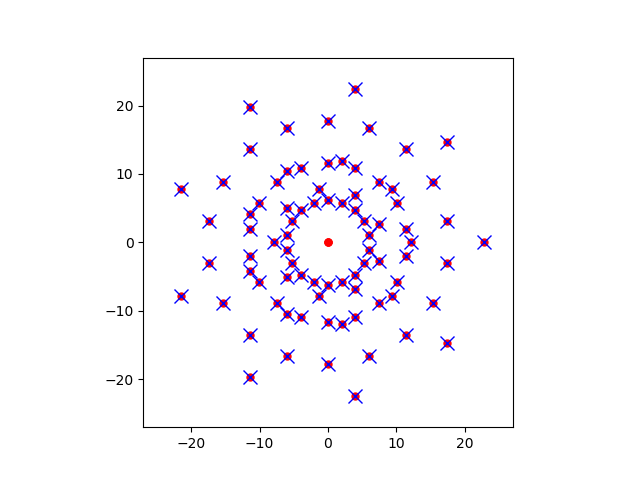}
        \caption{$\Gr(3,9)$}
        \label{fig:branes(3,9)}
    \end{subfigure}

    \caption{Summands of $\Fuk(\Gr(k,n))$ containing objects $T^{k(n-k)}_I$.}
    \label{fig:branes}
\end{figure}

In contrast with the case of projective spaces, for general Grassmannians there
are summands $\Fuk_\lambda(\Gr(k,n))$ with lower $|\lambda|$, and one can still use the
embedding $\theta_R$ above to find nonzero objects $T^{k(n-k)}_I$ in some of them, where
here $I$ is not a rotation of the special set of roots $I_0$; see Figure \ref{fig:branes}.
\par
One limitation is that one nonzero object might not be sufficient to generate
$\Fuk_\lambda(\Gr(k,n))$ when $|\lambda|$ is not maximum. A second limitation is that, depending on the arithmetic
of $k$ and $n$, the objects $T^{k(n-k)}_I$ can miss some summands of the Fukaya category. 
The simplest example of this phenomenon is $\Gr(2,4)$, for which the summand
$\Fuk_0(\Gr(2,4))$ of the Fukaya category does not contain any of the objects above.
Nohara-Ueda \cite{NU} showed in this case that the Gelfand-Cetlin system has
a monotone Lagrangian fiber whose topology is $S^1\times S^3$, and that it supports
generators for the $0$-summand. Earlier work of Smith \cite{Sm}
also implies that $\Fuk_0(\Gr(2,4))$ is generated by a vanishing cycle in a degeneration
of $\Gr(2,4)\subset\mathbb{P}^5$, thought as a quadric hypersurface.
\par
From our point of view, the objects $T^{k(n-k)}_I$ don't cover
all the summands of the Fukaya category whenever the critical points of $W$ are
not all contained in the rectangular torus chart of the mirror. In a separate work \cite{Cas}
we investigate how different charts of the mirror correspond to other monotone
Lagrangian tori in $\Gr(k,n)$.

Nevertheless Figure \ref{fig:branes} suggests a dichotomy: given a modulus,
none or many of the summands of the Fukaya category with that modulus contain objects $T^{k(n-k)}_I$
of the collection. We explain this phenomenon in terms of
\textit{equivariance} with respect to a suitable group action.

\begin{theorem}\label{thm:main-thm-2-statement}
If $T^{k(n-k)}_I$ is an object of $\Fuk_\lambda(\Gr(k,n))$, then it is nonzero and the objects
$$T^{k(n-k)}_{gI}\quad \text{in}\quad \Fuk_{g\lambda}(\Gr(k,n)) \quad \text{for all}\quad g\in D_n$$
are defined and nonzero as well.
\end{theorem}

For the proof of this theorem, most of the work goes into introducing an algebraic action of the dihedral group $D_n$ on the
Landau-Ginzburg mirror $\vGr(k,n)$ (Definitions \ref{DefDihedralRep} and \ref{DefYoungAction}) that makes the following commutative diagram
equivariant at the level of critical points and critical values (and globally $\mathbb{Z}/n\mathbb{Z}$-equivariant, where
$\mathbb{Z}/n\mathbb{Z}$ is the subgroup of $D_n$ generated by $r$):

\begin{equation}
\begin{tikzcd}[row sep=2cm, column sep=2cm]
^{D_n \circlearrowright}(\mathbb{C}^\times)^{k(n-k)} \arrow[rd, "W_{T^{k(n-k)}}" left] \arrow[r, hook, "\theta_R"] & \vGr(k,n)^{\circlearrowleft D_n}\arrow[d, "W"] \\
	& \mathbb{C}^{\circlearrowleft D_n}
\end{tikzcd}
\end{equation}

In the diagram, $W_{T^{k(n-k)}}$ denotes the Maslov 2 disk potential of the Gelfand-Cetlin
Lagrangian torus and commutativity is given by Theorem \ref{thm:main-thm-1-statement}. We call this action of $D_n$ on the Landau-Ginzburg mirror
\textit{Young action}, because it is defined in terms of Young diagrams.
In the Young action, the $s$ generator of $D_n$ does not act by conjugation:
indeed the action of $D_n$ is algebraic, whereas conjugation is not. On the other hand
the two actions do agree on the critical locus of the Landau-Ginzburg superpotential
$W$.
\par
The question of precisely what summands of the Fukaya category contain objects $T^{k(n-k)}_I$
appears to be related to the arithmetic of $k$ and $n$. When $n=p$ is prime, we give an argument
combining properties of \textit{vanishing sums} of roots of unity and Stanley's \textit{hook-content formula}
to show that one gets nonzero objects in all summands. In fact, in this case it also happens
that the quantum cohomology of $\Gr(k,p)$ has one dimensional summands, and this suffices to prove homological mirror symmetry.

\begin{theorem}\label{thm:main-thm-3-statement}
When $n=p$ is prime the objects $T^{k(p-k)}_I$ split-generate the Fukaya category of $\Gr(k,p)$, and
for every $\lambda\in\mathbb{C}$ there is an equivalence of triangulated categories
$$\D\Fuk_\lambda(\Gr(k,p))\simeq \D\mathcal{S}(W^{-1}(\lambda))\quad .$$
\end{theorem}

This article is the first step in a project to describe generators of the monotone Fukaya category
for all complex Grassmannians $\Gr(k,n)$. Other relevant works related to this problem are
the general approach to generation for Hamiltonian $G$-manifolds of Evans-Lekili \cite{EL} and
the study of immersed Lagrangians in Grassmannians of Hong-Kim-Lau \cite{HKL}.

\section*{Setup}\label{Setup}

A closed symplectic manifold $(X^{2N},\omega )$ is \textit{monotone} if $c_1(TX)=C[\omega]$
for some constant $C >0$, where the tangent bundle $TX\to X$ is $J$-complex with respect
to any almost complex structure $J$ on $X$ compatible with $\omega$. We will consider closed oriented
Lagrangian submanifolds $L\subset X$ that are themselves monotone, meaning that $\mu(D)=B\omega(D)$
for some constant $B >0$ and any $J$-holomorphic disk $D$ in $X$ with boundary on $L$, where
$\mu$ denotes the Maslov index. The constants $B$ and $C$ are related by $B=2C$.
\par
Working over an algebraically closed field $\mathbb{F}=\overline{\mathbb{F}}$, we denote by $\QH(X)$ the quantum
cohomology of $(X,\omega)$ over $\mathbb{F}$ (i.e. with Novikov parameter $q=1$). In the examples we are interested in the cohomology
of $X$ is concentrated in even degrees, so that $\QH(X)$ is a commutative unital
$\mathbb{F}$-algebra. The operator of quantum multiplication by the first Chern class $c_1\star$ induces a \textit{spectral decomposition}
$$\QH(X)=\bigoplus_{\lambda\in\mathbb{F}}\QH_\lambda(X)$$
in generalized eigenspaces, labelled by the eigenvalues. For each $\lambda$ one has a $\mathbb{Z}/2\mathbb{Z}$-graded 
$A_\infty$-category over $\mathbb{F}$ denoted $\Fuk_\lambda(X)$, where the objects $L_\xi$
are closed oriented monotone Lagrangian submanifolds $L\subset X$ equipped with a rank
one $\mathbb{F}$-linear local system $\xi$ whose holonomy $\operatorname{hol}_\xi:\pi_1(L)\to \mathbb{F}^\times$ satisfies
$$m^0(L_\xi)=\sum_{\beta}\#\mathcal{M}_J(L;\beta)\operatorname{hol}_\xi(\partial\beta )=\lambda\quad .$$
The sum above is over $\beta\in H_2(X,L;\mathbb{Z})$ with Maslov index $\mu(\beta)=2$, and $\#\mathcal{M}_J(L;\beta)$ denotes
the number of $J$-holomorphic disks through a generic point of $L$ in class $\beta$ for generic $\omega$-compatible $J$. This is a
finite sum thanks to the monotonicity assumptions. When $L\cong T^N$ is a torus, the fact that $\mathbb{F}^\times$ is abelian group allows us to think
$$\operatorname{hol}_\xi\in\operatorname{Hom}(H_1(T^N;\mathbb{Z}),\mathbb{F}^\times)\cong H^1(T^N;\mathbb{F}^\times)\cong(\mathbb{F}^\times)^N\quad .$$
The first isomorphism is natural, while the second depends on the choice of a basis $\gamma_1,\ldots ,\gamma_N$
of 1-cycles. For explicit calculations, it is convenient to specify such a basis and
call $x_1,\ldots ,x_N$ the relative coordinates, so that $\operatorname{hol}_\xi\mapsto m^0(L_\xi)$ gives an
algebraic function
$$W_{T^N}\in\mathcal{O}((\mathbb{F}^\times)^N)=\mathbb{F}[x_1^{\pm 1},\ldots , x_N^{\pm 1}]$$
that we call \textit{Maslov 2 disk potential} of the torus $T^N$. Disk potentials in different bases are related by
an integral linear change of variable, and the $GL(N,\mathbb{Z})$-orbit of $W_{T^N}$
is an Hamiltonian isotopy invariant of $T^N$.
\par
Given two objects $L^0_{\xi_0},L^1_{\xi_1}$ of $\Fuk_\lambda(X)$, the $\mathbb{F}$-vector space
$\hom(L^0_{\xi_0}, L^1_{\xi_1})$ has basis $L^0\cap L^1$ after an
Hamiltonian isotopy that arranges them to be transversal. The $\mathbb{Z}/2\mathbb{Z}$-grading
$$\hom(L^0_{\xi_0}, L^1_{\xi_1})=\hom^0(L^0_{\xi_0}, L^1_{\xi_1})\oplus\hom^1(L^0_{\xi_0}, L^1_{\xi_1})$$
is given by the 2-fold covering $\mathcal{L}^2(T_pX)\to\mathcal{L}(T_pX)$ of the
oriented Lagrangian Grassmannian of $T_pX$ over the unoriented one: the orientations on the
Lagrangians determine two points in the fibers over $T_pL^0$ and $T_pL^1$, with $|p|=0$ if
a \textit{canonical path} from $T_pL^0$ to $T_pL^1$ in $\mathcal{L}(T_pX)$ lifts
to a path connecting the orientations and $|p|=1$ otherwise. We adopt the convention that
the canonical path is given by $\pi/2$ counter-clockwise rotation in each $\langle \partial_{x_i},\partial_{y_i}\rangle$ plane in the local model
$$(T_pX,\omega_p)=(\mathbb{R}^{2N},dx_1\wedge dy_1+\ldots +dx_N\wedge dy_N)\quad ,$$
where $T_pL^0=\langle \partial_{x_1},\ldots ,\partial_{x_N}\rangle$ and $T_pL^1=\langle \partial_{y_1},\ldots ,\partial_{y_N}\rangle$.
For any $l\geq 1$ one has $\mathbb{F}$-linear maps
$$m^l:\hom(L^{l-1}_{\xi_{l-1}}, L^{l}_{\xi_{l}})\otimes\cdots\otimes\hom(L^0_{\xi_0}, L^1_{\xi_1})\to\hom(L^0_{\xi_0}, L^l_{\xi_l})$$
defined on intersection points by
$$m^l(p_l\otimes\cdots\otimes p_1)=\sum_{q\in L^0\cap L^l}\left( \sum_{\beta}\#\mathcal{M}_J(p_1,\ldots ,p_l,q;\beta)\operatorname{hol}_\xi(\partial\beta)\right)q\quad .$$
Here the sum is over $\beta\in H_2(X,L^0\cup \cdots \cup L^l;\mathbb{Z})$ whose Maslov index makes finite the count
$\#\mathcal{M}_J(p_1,\ldots ,p_l,q;\beta)$ of $J$-holomorphic strips with $p_1,\ldots ,p_l$ input and $q$ output asymptotic boundary conditions. Finally, 
$\operatorname{hol}_\xi(\partial\beta)\in\mathbb{F}^\times$ is obtained as compound holonomy
of the local systems $\xi_0,\ldots ,\xi_l$ along $\partial\beta$.
This construction depends on the choice of Hamiltonian isotopies and generic almost complex
structures $J$, but different choices produce equivalent $A_\infty$-categories and $\Fuk_\lambda(X)$ denotes
any of them. If $L^0_{\xi_0}$ and $L^1_{\xi_1}$ are monotone Lagrangians equipped with rank one
$\mathbb{F}$-linear local systems, we will always assume that
$m^0(L^0_{\xi_0})=m^0(L^1_{\xi_1})=\lambda$ and consider $\Fuk_\lambda(X)$ for different $\lambda\in\mathbb{F}$
as separate $A_\infty$-categories. This choice guarantees that $(m^1)^2=0$ in each morphism space,
so that we have well defined \textit{Floer cohomology} $\HF(L^0_{\xi_0},L^1_{\xi_1})$. This is in general only a vector
space over $\mathbb{F}$, but when $L^0_{\xi_0}=L^1_{\xi_1}$ it has an algebra structure induced
by the $A_\infty$-algebra structure of $\hom(L^0_{\xi_0},L^1_{\xi_1})$.
\par
We denote $\D\Fuk_\lambda(X)$ the \textit{derived category}. This is the homotopy category
of the enlargement of $\Fuk_\lambda(X)$ to the split-closure of the $A_\infty$-category of \textit{twisted complexes} $\operatorname{Tw}\Fuk_\lambda(X)$, where
sums of objects, shifts and cones of closed morphisms are available. The details of this construction
are not relevant here; it suffices to say that a set of objects of $\Fuk_\lambda(X)$ is
said to \textit{generate} $\D\Fuk_\lambda(X)$ whenever the smallest triangulated subcategory
containing the objects is the ambient category. For more details about the construction
of the monotone Fukaya category, we refer to Sheridan \cite{Sh}.
\par
In this article, we consider a Landau-Ginzburg mirror $X^\vee$ to $X$ which is a smooth
affine algebraic variety over $\mathbb{F}$, whose dimension $\operatorname{dim}(X^\vee)=N$ is half the
real dimension of the symplectic manifold $X$. Writing $X^\vee=\operatorname{Spec}(R)$ with
$R$ algebra over $\mathbb{F}$ of Krull dimension $N$, the superpotential $W:X^\vee\to\mathbb{F}$
will be an algebraic function $W\in R$. Smoothness of $X^\vee$ guarantees that the sheaf of
algebraic 1-forms $\Omega^1_{R/\mathbb{F}}$ is locally free of rank $N$, and the equation
$dW=0$ defines a closed subscheme $Z\subset X^\vee$, the \textit{critical locus},
which in our examples is always 0-dimensional. We call \textit{Jacobian ring} of $W$
the ring of algebraic functions on the critical locus $Z=\operatorname{Spec}(\operatorname{Jac}(W))$.
The fiber $W^{-1}(\lambda)$ over a closed point $\lambda\in\mathbb{F}$ is also a
closed subscheme, with $W^{-1}(\lambda)=\operatorname{Spec}(R/(W-\lambda))$.
The critical locus $Z$ decomposes as a union of 0-dimensional closed subschemes $Z_\lambda=Z\cap W^{-1}(\lambda)$,
and this induces a decomposition
$$\operatorname{Jac(W)}=\bigoplus_{\lambda\in\mathbb{F}}\operatorname{Jac}_\lambda(W)$$
that mirrors the spectral decomposition of the quantum cohomology $\QH(X)$, where we have
$Z_\lambda = \operatorname{Spec}(\operatorname{Jac}_\lambda(W))$ and $\operatorname{Jac}_\lambda(W)=\operatorname{Jac}(W)\otimes_{R} R/(W-\lambda)$.
\par
In this setting, for each $\lambda\in\mathbb{F}$ Orlov's derived category of singularities of $W^{-1}(\lambda)$
is equivalent to the homotopy category of the category of \textit{matrix factorizations} of $W-\lambda$
$$\D\mathcal{S}(W^{-1}(\lambda))\simeq \D\mathcal{M}(R,W-\lambda)\quad .$$
The category $\mathcal{M}(R,W-\lambda)$ is a \textit{differential graded category} with grading given by $\mathbb{Z}/2\mathbb{Z}$,
whose objects are $R$-modules $M=M_0\oplus M_1$ with finitely generated projective summands
of degree $0$ and $1$, and equipped with an $R$-linear map $d_M:M\to M$  of odd degree satisfying the equation
$(d_M)^2 = (W-\lambda)\operatorname{id}_M$. Morphisms between two matrix factorizations $M$ and $N$ in
$\mathcal{M}(R,W-\lambda)$ are given by cycles of a $\mathbb{Z}/2\mathbb{Z}$-graded complex over $\mathbb{F}$ with
$$\hom(M,N)=\hom^0(M,N)\oplus\hom^1(M,N)\quad \textrm{and} \quad d(f)=d_N\circ f - (-1)^{|f|}f\circ d_M \quad ,$$
where $f:M\to N$ denotes an $R$-linear map of degree $|f|$.
The name matrix factorizations comes from the fact that, if $M=M_0\oplus M_1$ with finitely
generated free summands, we can pick bases and represent $d_M$ as a matrix
\[
d_M=
\left(
\begin{array}{c|c}
0 & d_M^{01} \\
\hline
d_M^{10} & 0
\end{array}
\right)
\]
and the condition $(d_M)^2 = (W-\lambda)\operatorname{id}_M$ implies that $M_0$ and $M_1$
have same rank, with
$$d_M^{10}\circ d_M^{01} = (W-\lambda)\operatorname{id}_{M_0} \quad,\quad d_M^{01}\circ d_M^{10} = (W-\lambda)\operatorname{id}_{M_1}\quad .$$
Essentially, matrix factorizations encode the fact that any coherent sheaf on $W^{-1}(\lambda)$
admits a locally free resolution that eventually becomes 2-periodic, as proved by Eisenbud \cite{Ei}.
For more details about the category of matrix factorizations, we refer to Dyckerhoff \cite{Dy}.
\par
We list below for future reference some facts that will be crucial for this article.

\begin{theorem}\label{Fact1}(Auroux \cite[Proposition 6.8]{AuT}) If $\HF(L_\xi, L_\xi)\neq 0$, then $m^0(L_\xi)$ is
an eigenvalue of the operator $c_1\star$ acting on $\QH(X)$.
\end{theorem}

This tells us that the only $\lambda\in\mathbb{F}$ for which $\D\Fuk_\lambda(X)$ can be nontrivial
are those appearing in the spectral decomposition of quantum cohomology. 

\begin{theorem}\label{Fact2}(Auroux \cite[Proposition 6.9]{AuT}, Sheridan \cite[Proposition 4.2]{Sh}) If $T^N$ is a monotone Lagrangian torus, then
$\operatorname{hol}_\xi$ is a critical point of $W_{T^N}$ if and only if $\HF(T^N_\xi,T^N_\xi)\neq 0$.
\end{theorem}

This reduces the problem of showing that $T^N_\xi$ is a nonzero object of the Fukaya category to studying the critical
points of the Maslov 2 disk potential $W_{T^N}$.

\begin{theorem}\label{Fact3}(Sheridan \cite[Corollary 2.19]{Sh}) If $\QH_\lambda(X)$ is one-dimensional, any
object $L_\xi$ of $\Fuk_\lambda(X)$ with $\HF(L_\xi,L_\xi)\neq 0$ split-generates $\D\Fuk_\lambda(X)$.
\end{theorem}

The \textit{generation criterion} above is an adaptation of the
one introduced by Abouzaid \cite{Ab} to the monotone setup.

\begin{theorem}\label{Fact4}(Sheridan \cite[Proposition 4.2-4.3 and Corollary 6.5]{Sh})
If $\mathbb{F}=\mathbb{C}$ and $T^N$ is a monotone Lagrangian torus, then for every nondegenerate critical point $\operatorname{hol}_\xi$ of $W_{T^N}$
$$\HF(T^N_\xi,T^N_\xi)\cong \Cl_N$$
as $\mathbb{C}$-algebras, where $\Cl_N$ denotes the Clifford algebra of the quadratic form
of rank $N$ on $\mathbb{C}^N$. Moreover if $T^N_\xi$ generates $\Fuk_\lambda(X)$ there
is an equivalence of triangulated categories
$$\D\Fuk_\lambda(X)\simeq \D(\Cl_N)$$
with the derived category of finitely generated modules over $\Cl_N$.
\end{theorem}

We will use the theorem above as an ingredient to establish homological mirror symmetry
for Grassmannians $\Gr(k,p)$ with $p$ prime.

\begin{remark}
Working in the monotone setting allows us to ignore the \textit{Novikov field} $\Lambda_\mathbb{F}$
and use instead $\mathbb{F}$ directly, by setting the Novikov parameter $q=1$. More specifically, in this
paper we will only consider $\mathbb{F}=\mathbb{C}$.
\end{remark}

\textbf{Acknowledgements} I thank my PhD advisor Chris Woodward for suggesting to think about
Grassmannians and useful conversations. I also thank: Mohammed Abouzaid for helpful conversations
about mirror symmetry and the Fukaya category; Lev Borisov for clarifications about toric
degenerations; Anders Buch for showing me the quantum Pieri rule; Hiroshi Iritani for a
useful conversation about part (3) of Propositon \ref{Prop1}; Alex Kontorovich and
Stephen Miller for remarks on vanishing sums of roots of unity; Greg Moore for pointing
out relevant physics literature. I thank the referees for valuable comments, and
for suggesting a shorter proof of Lemma \ref{lemma:dihedral-action}.

\section{Closed mirror symmetry for $\Gr(k,n)$}\label{Section1}

The quantum cohomology of $\Gr(k,n)$ was first computed by Bertram \cite{Be}, and admits
a purely combinatorial description. As vector space, it agrees with the classical
cohomology $H^\bullet(\Gr(k,n);\mathbb{C})$. This has a basis given by Poincar\'{e}
dual classes $\sigma_d$ to certain algebraic cycles $X_d\subseteq \Gr(k,n)$ known
as \textit{Schubert varieties}; see for example Fulton \cite{Fu}. The classes
$\sigma_d$ are parametrized by \emph{Young diagrams} $d$ in a $k\times (n-k)$ grid,
simply denoted $d\subseteq k\times (n-k)$.

\begin{example}

If $k=2$ and $n=5$, the diagrams $d\subseteq 2\times 3$ are:

\ytableausetup{boxsize=1em}
\begin{center}
\vspace{0.3cm}
$\emptyset$ \; , \;
\ydiagram{1,0} \; , \;
\ydiagram{2,0} \; , \;
\ydiagram{1,1} \; , \;
\ydiagram{3,0} \; , \;
\ydiagram{2,1} \; , \;
\ydiagram{3,1} \; , \;
\ydiagram{2,2} \; , \;
\ydiagram{3,2} \; , \;
\ydiagram{3,3} \; .
\vspace{0.3cm}
\end{center}

\end{example}

If the $i$-th row has $d_i$ boxes, the condition for being a diagram is that $n-k\geq d_1\geq \ldots \geq d_k\geq 0$.
The tuple $d=(d_1,\ldots ,d_k)$ is a partition of the number of boxes $|d|$, which equals the
complex codimension of the corresponding Schubert variety $X_d$. It is also useful to think $d$ as
a subset of $\{1,\ldots ,n\}=[n]$, in which case $d^|$ denotes the size $k$ set of the vertical
steps, and $d^-$ the size $n-k$ set of horizontal steps. The steps of a diagram are
counted from the top-right corner of the grid (not of the diagram) to the bottom-left, following
the part of the border of $d$ in the interior of the grid.

\begin{example}

Let $k=3$ and $n=8$. Then $d=(3,2,0)$ denotes a Young diagram, thought as partition
of $|d|=5=3+2+0$. The vertical steps are $d^|=\{3,5,8\}$, and the horizontal steps
are $d^-=\{1,2,4,6,7\}$. These sets are computed by placing $d$ in the $3\times 5$
grid as follows:

	\vspace{0.3cm}
	\centering
	\includegraphics[width=0.25\textwidth]{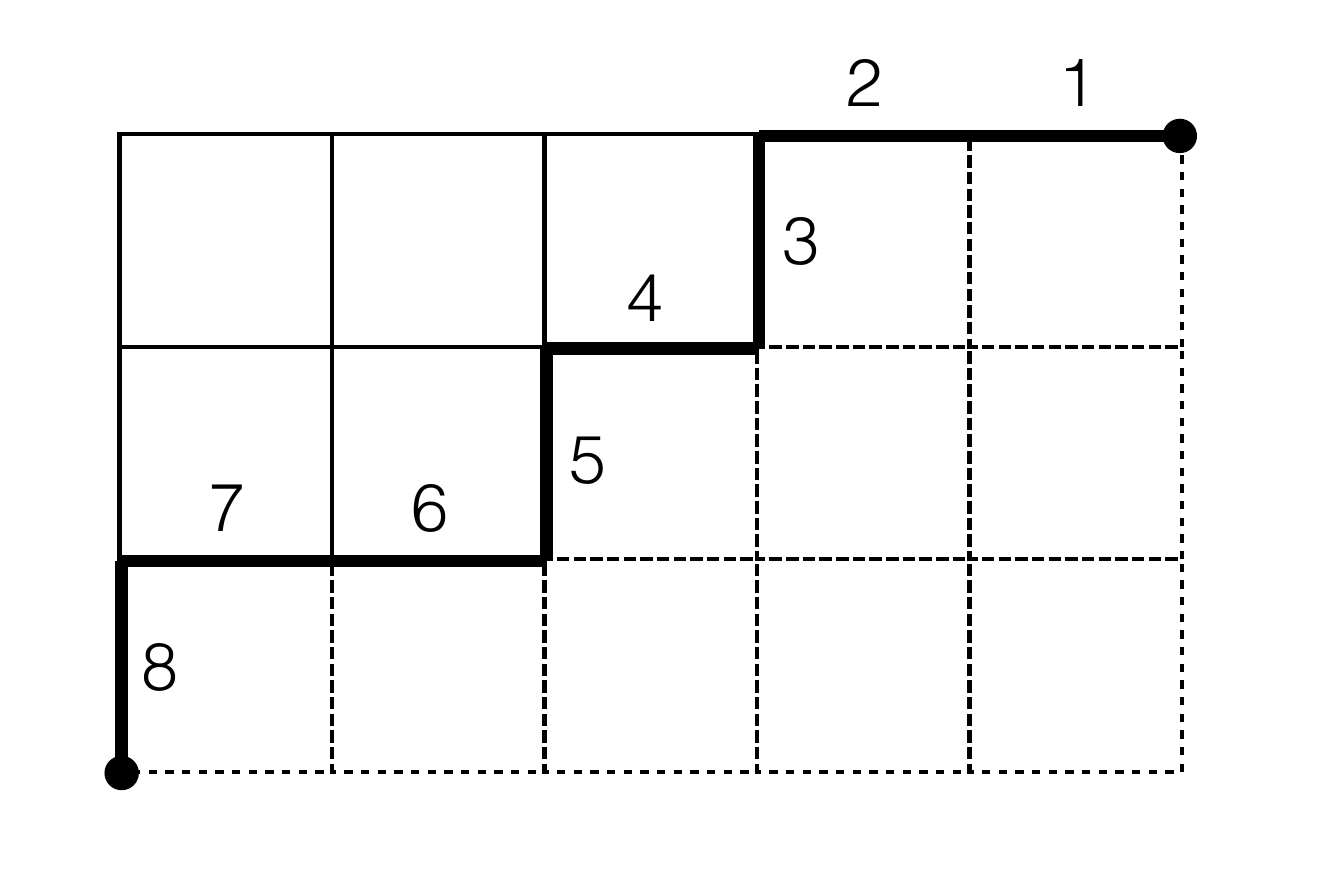}
	\vspace{0.3cm}
	
\end{example}

\ytableausetup{boxsize=0.2em}
We will not need the full ring structure of $\QH(\Gr(k,n)$. If $\star$ denotes
its product, it will suffice to mention that the product by $c_1(T\Gr(k,n))=n\sigma_{\ydiagram{1}}$,
is determined by the \textit{quantum Pieri rule}:
$$\sigma_{\ydiagram{1}}\star\sigma_d = \sigma_{\ydiagram{1}}\cdot\sigma_d + \sigma_{\hat{d}} \quad .$$
The first term in this formula is the cup product, while the second is a quantum correction.
The classical part $\sigma_{\ydiagram{1}}\cdot\sigma_d$ is a sum of Schubert classes
obtained by adding one box to $d$ in all possible ways. The quantum part $\sigma_{\hat{d}}$ is
a single Schubert class, with diagram $\hat{d}$ obtained by ereasing the
full first row and the full first colum of $d$, or $0$ otherwise (i.e. if $d$ has less than $n-k$ boxes
in first row or less than $k$ boxes in first column).

\begin{example}

Let $k=2$ and $n=5$. Some sample classical/quantum contributions are:
\ytableausetup{boxsize=1em}
\vspace{0.3cm}
$$\ydiagram{1,0} \cdot \ydiagram{2,1} = \ydiagram{3,1} + \ydiagram{2,2}\quad ,\quad \widehat{\ydiagram{3,2}} = \ydiagram{1,0}
\quad , \quad \widehat{\ydiagram{2,2}} = 0 \quad .$$
\vspace{0.3cm}

\end{example}

Building on Peterson's work on the quantum cohomology of flag varieties, Marsh-Rietsch \cite{MR} proposed
for $X=\Gr(k,n)$ a Landau-Ginzburg mirror which is an open subvariety of the dual Grassmannian
$X^\vee=\vGr(k,n)\subset\Gr(n-k,n)$, the complement of an explicit divisor
with a suitable $W:\vGr(k,n)\to\mathbb{C}$ given as rational function on $\Gr(n-k,n)$ with
poles on the divisor. Before getting to the details, we recall the basic facts about Pl\"{u}cker
coordinates.
\par
Think of $[M]\in\Gr(k,n)$ as full rank $k\times n$ matrix $M$ with complex entries, with
rows giving a basis of a $k$-dimensional linear subspace of $\mathbb{C}^n$; the basis is not unique, and
$[M]$ denotes the equivalence class of $M$ modulo row operations. For each Young diagram $d\subseteq k\times (n-k)$,
the \textit{Pl\"{u}cker coordinate} $p_d(M)$ denotes the determinant
of the minor of $M$ obtained by selecting the columns $d^|\subset [n]$. Ordering the Young diagrams $d$
according to the lexicographic order on the sets of vertical steps $d^|$, we get a map
$$\Gr(k,n)\to\mathbb{P}^{{n \choose k}-1}\quad , \quad [M] \mapsto [p_d(M)]$$
which is well defined because row operations on $M$ change all Pl\"{u}cker coordinates by the same
nonzero factor, and furthermore at least one Pl\"{u}cker coordinate is nonzero because $M$ is full rank.
The map above is a closed embedding of schemes, and there is also a \textit{dual} embedding
$$\Gr(n-k,n)\to\mathbb{P}^{{n \choose k}-1}\quad , \quad [\check{M}] \mapsto [\check{p}_d(\check{M})]\quad .$$
Here $[\check{M}]\in\Gr(n-k,n)$ denotes a full rank $n\times (n-k)$ matrix $\check{M}$ with complex entries modulo
column operations. Our convention is to always use diagrams $d\subseteq k\times (n-k)$, but to think
the dual Pl\"{u}cker coordinates $\check{p}_d(\check{M})$ as given by the determinant of the minor of $\check{M}$
obtained by selecting the rows $d^-\subseteq [n]$, ordering the diagrams according to the lexicographic
order on the sets of horizontal steps $d^-$.
\par
The Landau-Ginzburg mirror of $\Gr(k,n)$ is:
\ytableausetup{boxsize=0.2em}
$$\vGr(k,n)=\Gr(n-k,n)\setminus\{\check{p}_1\cdots\check{p}_n=0\}\quad , \quad
W = \frac{\check{p}_1^{\ydiagram{1}}}{\check{p}_1} + \ldots + \frac{\check{p}_n^{\ydiagram{1}}}{\check{p}_n}\quad .$$
Here $\check{p}_1, \ldots , \check{p}_n$ denote the Pl\"{u}cker coordinates of the $n$ \textit{boundary rectangular}
Young diagrams, where the $i$-th of them is called $d_i$ and has horizontal steps given by taking $i$ cyclic shifts
of the set $\{1,\ldots ,n-k\}$. We abbreviate $\check{p}_{d_i}=\check{p}_i$.
Observe that there are rectangular Young diagrams that are not boundary rectangular,
and we call them \textit{interior rectangular}. They are characterized by the fact that
none of the edges has full length, i.e. length $k$ or $n-k$;

\begin{example}

The boundary rectangular diagrams for $k=2$ and $n=5$ are:

\ytableausetup{boxsize=1em}
\vspace{0.3cm}
$$d_1=\emptyset \quad , \quad d_2=\ydiagram{3,0} \quad , \quad d_3=\ydiagram{3,3} \quad , \quad d_4=\ydiagram{2,2} \quad , \quad d_5=\ydiagram{1,1} \quad ;$$

\vspace{0.3cm}
their horizontal steps are $d_1^-=123$, $d_2^-=234$, $d_3^-=345$, $d_4^-=145$, $d_5^-=125$.

\end{example}

\begin{example}

The interior rectangular diagrams for $k=3$ and $n=7$ are:

\vspace{0.5cm}
$$\ydiagram{1,0,0} \quad , \quad \ydiagram{2,0,0} \quad , \quad \ydiagram{3,0,0} \quad , \quad \ydiagram{1,1,0} \quad , \quad \ydiagram{2,2,0} \quad , \quad \ydiagram{3,3,0} \quad .$$

\vspace{0.5cm}
Note that $\emptyset$ is considered boundary rectangular by convention.

\end{example}

\ytableausetup{boxsize=0.2em}
Going back to the formula for the superpotential $W$, $\check{p}_i^{\ydiagram{1}}$ denotes the Pl\"{u}cker coordinate of the Young diagram obtained
by using the quantum Pieri rule to compute the $\star$ product of
\ytableausetup{boxsize=0.5em}
$\ydiagram{1}$ with the $i$-th boundary rectangular diagram. Note that for these the product
with $\ydiagram{1}$ is a single diagram instead of a sum of diagrams.
\ytableausetup{boxsize=0.2em}

\begin{example}

The Landau-Ginzburg mirror of $\Gr(2,5)$ is given by
$$\vGr(2,5)=\Gr(3,5)\;\setminus\;\{ \; \check{p}_{\emptyset}\check{p}_{\ydiagram{3,0}}\check{p}_{\ydiagram{3,3}}\check{p}_{\ydiagram{2,2}}\check{p}_{\ydiagram{1,1}}=0 \; \} \quad ,$$
with superpotential
$$W\quad =\quad \frac{\check{p}_{\ydiagram{1,0}\star\emptyset}}{\check{p}_{\emptyset}} +
\frac{\check{p}_{\ydiagram{1,0}\star\ydiagram{3,0}}}{\check{p}_{\ydiagram{3,0}}} +
\frac{\check{p}_{\ydiagram{1,0}\star\ydiagram{3,3}}}{\check{p}_{\ydiagram{3,3}}} +
\frac{\check{p}_{\ydiagram{1,0}\star\ydiagram{2,2}}}{\check{p}_{\ydiagram{2,2}}} +
\frac{\check{p}_{\ydiagram{1,0}\star\ydiagram{1,1}}}{\check{p}_{\ydiagram{1,1}}}\quad = \quad
\frac{\check{p}_{\ydiagram{1,0}}}{\check{p}_{\emptyset}} +
\frac{\check{p}_{\ydiagram{3,1}}}{\check{p}_{\ydiagram{3,0}}} +
\frac{\check{p}_{\ydiagram{2,0}}}{\check{p}_{\ydiagram{3,3}}} +
\frac{\check{p}_{\ydiagram{3,2}}}{\check{p}_{\ydiagram{2,2}}} +
\frac{\check{p}_{\ydiagram{2,1}}}{\check{p}_{\ydiagram{1,1}}} \quad .$$

\end{example}

\begin{remark}

The Landau-Ginzburg mirror $\vGr(k,n)=\Gr(n-k,n)\setminus\{\check{p}_1\cdots\check{p}_n=0\}$ is an affine variety
because it is the complement of an ample divisor in a projective variety. We will call its
coordinate ring $\mathbb{C}[\vGr(k,n)]$.

\end{remark}

The following result says that this is a correct mirror for $\Gr(k,n)$ in
the sense of closed mirror symmetry.

\begin{theorem}\label{ClosedMS}(Rietsch \cite[Theorem 4.1]{Ri2})
There is an isomorphism of $\mathbb{C}$-algebras
$$\QH(\Gr(k,n))\cong \operatorname{Jac}(W)$$
where the Jacobian ring of $W$ on the right is given by
$\operatorname{Jac}(W)=\mathbb{C}[\vGr(k,n)]/(\partial_{\check{p}_d}W,\forall d)$.
\end{theorem}

The final goal of this section is to determine the spectral decomposition of $\QH(\Gr(k,n))$ labelled
by eigenvalues of the operator of quantum multiplication by the first Chern class $c_1\star$.
The eigenvalues with maxmimum modulus in the spectrum of this operator
have been studied before. In fact, for general smooth projective Fano varieties,
they are expected to satisfy the following conjecture.

\begin{conjecture}(Galkin-Golyshev-Iritani \cite[Conjecture 3.1.2]{GGI})
If $X$ is a smooth projective Fano variety, denote
$$T = \max \{ \; |\lambda | \; : \; \lambda \; \textrm{eigenvalue of} \; c_1\star \; \}$$
and $r$ its Fano index. Then:
\begin{enumerate}
	\item T is an eigenvalue of $c_1\star$ ;
	\item if $\lambda$ is an eigenvalue of $c_1\star$ with $|\lambda|=T$, then
		$u=T\zeta$ for some $r$-th root of unity $\zeta\in\mathbb{C}$ ;
	\item the multiplicity of $T$ is 1 .
\end{enumerate}
\end{conjecture}

The conjecture above is known to hold for Grassmannians and more general
homogeneous varieties $G/P$; see Galkin-Golyshev-Iritani \cite{GGI} and Cheong-Li \cite{CL}.
Since we are interested in the full spectrum of $c_1\star$, it is in principle possible to
use the quantum Pieri rule presented earlier to compute it in specific cases.
However, in the proof of Proposition \ref{Prop1} we use a different approach, relying on
the existence of a particular basis for $\QH(\Gr(k,n))$: the \textit{Schur basis}. The
author learned of this basis from the work of Rietsch \cite{Ri}.
\par
To each $d\subseteq k\times (n-k)$ one can associate a symmetric polynomial
in $k$ variables, called \textit{Schur polynomial} of $d$, and defined as
$$S_d(x_1,\ldots ,x_k) = \sum_{T_d}x_1^{t_1}\cdots x_k^{t_k}\quad .$$
The sum is over \emph{semi-standard Young tableau} on the diagram $d$, obtained by filling $d$
with labels $\{1,\ldots ,k\}$ in such a way that rows are weakly increasing and columns are strictly
increasing. The exponent $t_i$ of $x_i$ records the number of occurrences of the label $i$.

\begin{example}

The following is an example with $k=2$:
\ytableausetup{boxsize=normal}
$$ d = \ydiagram{3,1}\quad ; \quad T_d =
\begin{ytableau}
 1 & 1 & 1 \\
 2 \\
\end{ytableau} \quad ,\quad
\begin{ytableau}
 1 & 1 & 2 \\
 2 \\
\end{ytableau} \quad ,\quad
\begin{ytableau}
 1 & 2 & 2 \\
 2 \\
\end{ytableau}\quad .$$
The corresponding Schur polynomial is
\ytableausetup{boxsize=0.2em}
$S_{\ydiagram{3,1}}(x_1,x_2)=x_1^3x_2 + x_1^2x_2^2 + x_1x_2^3$ .

\end{example}

The Schur basis $\sigma_I$ 
is indexed by sets $I$ with $|I|=k$ of roots of $x^n = (-1)^{k+1}$ and given by
$$\sigma_I =\sum_{d}\overline{S_d(I)}\sigma_d\quad .$$
The elements of the Schur basis are eigenvectors for $c_1\star$,  and in fact for
any operator $\sigma_d\star$; see \cite{Ri}. For the reader's convenience, we
give a proof below.

\begin{lemma}(Rietsch \cite{Ri})
\label{lemma:schur-basis}
For any $d\subseteq k\times (n-k)$ and any set $I$ of distinct
roots of $x^n=(-1)^{k+1}$, the following identity holds:
$$\sigma_d\star \sigma_I = S_d(I)\sigma_I \quad .$$
\end{lemma}

\begin{proof}
By definition of Schur basis
$$\sigma_d\star\sigma_I = \sum_\nu\overline{S_\nu(I)}\sigma_d\star\sigma_\nu \quad ,$$
where $\nu\subseteq k\times (n-k)$ runs among all Young diagrams. The product
of two Schubert classes in $\QH(\Gr(k,n))$ is given by
$$\sigma_d\star\sigma_\nu = \sum_\mu\left( \sum_{h\in\mathbb{N}}\langle \sigma_d, \sigma_\nu, \sigma_{PD(\mu)}\rangle_h q^h\right)\sigma_\mu \quad ,$$
where $PD(\mu)$ denotes the Poincar\'{e} dual Young diagram of $\mu$, and
$\langle \sigma_d, \sigma_\nu, \sigma_{PD(\mu)}\rangle_h$ is the
Gromov-Witten invariant counting genus $0$ curves of degree $h$ through the Schubert
cycles $X_d, X_\nu, X_{PD(\mu)}\subseteq \Gr(k,n)$. This number is $0$ unless
$|d|+|\nu|+|PD(\mu)| = k(n-k) + hn$, so that calling
$$h(\mu) = \frac{|d|+|\nu|+|PD(\mu)|-k(n-k)}{n}$$
we get
$$\sigma_d\star\sigma_\nu = \sum_\mu \langle \sigma_d, \sigma_\nu, \sigma_{PD(\mu)} \rangle_{h(\mu)} q^{h(\mu)}\sigma_\mu \quad .$$
The Gromov-Witten invariant above has a very explicit expression, known as
Bertram-Vafa-Intriligator formula:
$$\langle \sigma_d, \sigma_\nu, \sigma_{PD(\mu)} \rangle_{h(\mu)} = \frac{1}{n^k}\sum_{J}\frac{S_d(J)S_\nu(J)S_{PD(\mu)}(J)}{S_{k\times (n-k)}(J)}\vert\operatorname{Van}(J)\vert^2 \quad ;$$
see \cite[Theorem 10.3]{Ri} (also \cite{Be}, \cite{ST}). In the formula, the sum
runs over $J$ size $k$ sets of roots of $x^n=(-1)^{k+1}$, and if $J=\{\zeta_1, \ldots ,\zeta_k\}$
$$\vert\operatorname{Van}(J)\vert = \prod_{i<j}|\zeta_i - \zeta_j| \quad .$$
Rearranging the sums in the product $\sigma_d\star\sigma_I$, and using the fact
(\cite[Lemma 4.4]{Ri}) that
$$\frac{S_{PD(\mu)}(J)}{S_{k\times (n-k)}(J)} = \overline{S_\mu(J)} \quad ,$$
we find that
$$\sigma_d\star\sigma_I = \frac{1}{n^k}\sum_{\mu}q^{h(\mu)}\sigma_\mu\left( \sum_{J} S_d(J)\overline{S_\mu(J)}
\vert\operatorname{Van}(J)\vert^2 \left( \sum_{\nu} S_\nu(J)\overline{S_\nu(I)}\right) \right) \quad .$$
Now \cite[Proposition 4.3]{Ri} says that
$$\sum_\nu S_\nu(J)\overline{S_\nu(I)} = \frac{n^k}{\vert\operatorname{Van}(J)\vert^2}\delta_{J,I} \quad ,$$
so we arrive to
$$\sigma_d\star\sigma_I = S_d(I) \left( \sum_\mu q^{h(\mu)}\overline{S_\mu(I)}\sigma_\mu \right) \quad .$$
Evaluating at $q=1$ we find the desired formula.
\end{proof}

We are now ready to prove  the following.

\begin{proposition}\label{Prop1}
The following properties hold:
\begin{enumerate}
	\item The eigenvalues of $c_1\star$ acting on $\QH(\Gr(k,n))$ are given by $n(\zeta_1+\ldots +\zeta_k)$, 
	with $\{\zeta_1,\ldots \zeta_k\}$ varying among the size $k$ sets of roots of $x^n = (-1)^{k+1}$.
	\item Let $O(2)$ act on the complex plane by linear isometries of the Euclidean metric.
		The subgroup that maps the set of eigenvalues of $c_1\star $ to itself is isomorphic to the dihedral group $D_n$.
	\item If $n=p$ prime, then all eigenvalues of $c_1\star$ have multiplicity one.
\end{enumerate}
\end{proposition}

\begin{proof}

\textit{1)} Follows from $c_1=n\sigma_{\ydiagram{1}}$ and the fact that a single box Young diagram supports exactly
$k$ tableaux, obtained by labelling it with any of the labels in $\{1,\ldots ,k\}$, so that
$$S_{\ydiagram{1}}(x_1,\ldots ,x_k)=x_1+\ldots +x_k\quad .$$
\textit{2)} If $I=\{\zeta_1,\ldots ,\zeta_k\}$, rotation of $2\pi /n$ and conjugation give
$$e^{2\pi i/n}nS_{\ydiagram{1}}(I)=n(e^{2\pi i/n}\zeta_1+\ldots +e^{2\pi i/n}\zeta_k)=nS_{\ydiagram{1}}(e^{2\pi i/n}I)$$
$$\overline{nS_{\ydiagram{1}}(I)}=n(\overline{\zeta_1}+\ldots +\overline{\zeta_k})=nS_{\ydiagram{1}}(\overline{I})$$
and these two trasformations generate a copy of $D_n$ in the subgroup of $O(2)$ that preserves
the eigenvalues. There are no other transformations with this property because the subgroup is finite, and
the only finite subgroups of $O(2)$ are cyclic or dihedral; therefore it must be contained
in a dihedral group, possibly larger than $D_n$. On the other hand, it cannot be larger than $D_n$
because there are $n$ eigenvalues with maximum modulus: this follows from the fact that $n$ is the Fano index of $\Gr(k,n)$, and
the Grassmannians have property $\mathcal{O}$ introduced by Galkin-Golyshev-Iritani \cite{GGI} (see also \cite[Proposition 3.3 and Corollary 4.11]{CL}), so that any element in our subgroup must be in
particular a symmetry of the $n$-gon formed by the eigenvalues of maximum modulus.\\\\
\textit{3)} For $p=2$ we must have $k=1$, and the statement is true. If $p>2$ prime,
	observe that the statement for $\Gr(k,p)$ is equivalent to the one for $\Gr(p-k,p)$ because the
	two Grassmannians are isomorphic. We can use this to assume without loss of generality that
	$k$ is odd, since when it is even we can replace $k$ with $p-k$. Let now $\{\xi_1,\ldots ,\xi_k\}$ and $\{\zeta_1,\ldots ,\zeta_k\}$ be two size $k$
	sets of $p$-th roots of
	$$x^p = (-1)^{k+1} = 1$$
	and call $z=e^{2\pi i/p}$. Rewrite
	$$\xi_1+\ldots +\xi_k = z^{i_1} + \ldots + z^{i_k}\quad ,\quad \zeta_1+\ldots +\zeta_k = z^{j_1} + \ldots + z^{j_k}$$
	with $0\leq i_1<\ldots <i_k <p$ and $0\leq j_1<\ldots <j_k <p$. Denote $\langle z \rangle$ the subgroup of $\mathbb{C}^\times$
	generated by $z$. The map $\phi(1)=z$
	extends to a morphism of group rings
	$$\phi : \mathbb{Z}[\mathbb{Z}/p\mathbb{Z}]\to\mathbb{Z}[\langle z\rangle ]\quad .$$
	Think now of the sums above as elements of a group ring
	$$z^{i_1} + \ldots + z^{i_k} = \phi\left(\sum_{u\in\mathbb{Z}/p\mathbb{Z}}a_ut^u\right)=\phi(a)\quad ,\quad 
	z^{j_1} + \ldots + z^{j_k} = \phi\left(\sum_{u\in\mathbb{Z}/p\mathbb{Z}}b_ut^u\right)=\phi(b)$$
	where $a,b\in\mathbb{Z}[\mathbb{Z}/p\mathbb{Z}]$ have coefficients
	$$	a_u = \begin{cases}
				1 &\text{if } u\in\{i_1,\ldots ,i_k\}\\
				0 &\text{ otherwise}
				\end{cases}	\quad , \quad
		b_u = \begin{cases}
				1 &\text{if } u\in\{j_1,\ldots ,j_k\}\\
				0 &\text{ otherwise}
				\end{cases}	\quad\quad .
	$$
	Now the two eigenvalues of $c_1\star $ corresponding to $\{i_1,\ldots ,i_k\}$ and $\{j_1,\ldots ,j_k\}$
	are equal whenever $\phi(a)=\phi(b)$, or equivalently $\phi(a-b)=0$. The kernel of
	the morphism $\phi$ is known (de Bruijn \cite[Theorem 1]{dB}; see also Lam-Leung \cite[Theorem 2.2]{LL}) and it is
	$$\operatorname{ker}\phi = \{ l(1+t+\ldots +t^{p-1})\, :\, l\in\mathbb{Z}\}\quad .$$
	Therefore there exists $l\in\mathbb{Z}$ such that
	$$\sum_{u\in\mathbb{Z}/p\mathbb{Z}}(a_u-b_u)t^u = l + lt + \ldots + lt^{p-1}$$
	so that $a_u-b_u=l$ for every $u\in\mathbb{Z}/p\mathbb{Z}$. Observe that for every $u$
	we have $a_u-b_u\in\{-1,0,1\}$, and moreover we can't have $l=\pm 1$ because both $a$ and
	$b$ have exactly $k$ of their coefficients ($a_u$ and $b_u$ respectively) different from $0$.
	We conclude that $l=0$, so that $a=b$ and therefore $\{i_1,\ldots ,i_k\}=\{j_1,\ldots ,j_k\}$.
\end{proof}

\section{Critical points and torus charts}\label{Section2}

The critical points of the Landau-Ginzburg superpotential $W:\vGr(k,n)\to\mathbb{C}$
have been computed by Rietsch.
\begin{theorem}\label{ThmKa}(\cite[Proposition 9.3]{MR}, \cite[Theorem 3.4 and Lemma 3.7]{Ri}; see also \cite[Theorem 1.1 and Corollary 3.12]{Ka})
The critical points of $W:\vGr(k,n)\to\mathbb{C}$ are given by
$$[\check{M}_{I^\vee}]=
\begin{bmatrix}
    1 & 1 & 1 & \dots  & 1 \\
    \zeta_1 & \zeta_2 & \zeta_3 & \dots  & \zeta_{n-k} \\
    \zeta_1^2 & \zeta_2^2 & \zeta_3^2 & \dots  & \zeta_{n-k}^2 \\
    \vdots & \vdots & \vdots &  & \vdots \\
    \zeta_1^{n-1} & \zeta_2^{n-1} & \zeta_3^{n-1} & \dots  & \zeta_{n-k}^{n-1}
\end{bmatrix}$$
where $I^\vee=\{\zeta_1,\ldots ,\zeta_{n-k}\}$ is a set of $n-k$ distinct roots of $x^n = (-1)^{n-k+1}$.
\end{theorem}

Observe that the matrix above is always full rank because, the roots being distinct,
by the Vandermonde formula we have
$$\check{p}_1([\check{M}_{I^\vee}])=\check{p}_{\emptyset}([\check{M}_{I^\vee}])=\prod_{1\leq i<j\leq n-k}(\zeta_j-\zeta_i)\neq 0\quad .$$

According to Scott \cite{Sc} (see also Rietsch-Williams \cite{RW}), $\vGr(k,n)$
contains open subschemes that are algebraic tori of the form
$$T_\cC=\{\, [\check{M}]\in\vGr(k,n)\, :\, \check{p}_d([\check{M}])\neq 0\, \forall d\in\cC\, \}$$
labeled by a collection of Young diagrams $\cC$
in a $k\times (n-k)$ grid. Each collection $\cC$ is such that the Pl\"{u}cker coordinates
$\check{p}_d$ with $d\in\cC$ form a transcendence basis of the function field
$\mathbb{C}(\vGr(k,n))$, and these can be thought as coordinates for the
torus chart $T_\cC\cong (\mathbb{C}^\times)^{k(n-k)}$.
Restricting the superpotential $W$ to each torus chart, we get Laurent polynomials
$$W_\cC=W\restriction_{T_\cC}\in\mathbb{C}[\,\check{p}_d^{\pm 1}\,: \, d\in\cC \,]\quad .$$
If $\cC$ and $\cC'$ are different collections of Pl\"{u}cker coordinates, the corresponding Laurent polynomials
$W_{\cC}$ and $W_{\cC'}$ involve different sets of variables $\check{p}_d$: for $W_\cC$ the variables are Pl\"{u}cker coordinates
parametrized by Young diagrams $d\in\cC$, while for $W_{\cC'}$ by diagrams $d\in\cC'$.
\par
For the purposes of this paper, we will focus only on one torus chart: the \textit{rectangular chart} $T_{\cC^R}$.
See \cite{Cas} for a symplectic topology interpretation of the other charts.
The name comes from the fact that it corresponds to the collection $\cC^R$
of rectangular Young diagrams in the $k\times(n-k)$ grid (excluding $\emptyset$, which
can be thought as corresponding to setting the corresponding Pl\"{u}cker coordinate
$\check{p}_\emptyset = 1$).
\par
This chart has the property that it always exists in $\vGr(k,n)$, no matter what $k$ and $n$ are; instead an arbitrary collection
of $k(n-k)$ Young diagrams in a $k\times(n-k)$ grid does not in general give a
tanscendence basis for $\mathbb{C}(\vGr(k,n))$. One more consideration about the rectangular torus
chart $T_{\cC^R}$ is that when $k=1$, i.e. for $\Gr(1,n)=\mathbb{P}^{n-1}$,
$\vGr(1,n)=(\mathbb{C}^\star)^{n-1}$ is a single torus, and there
is only one possible way to choose $1\cdot(n-1)=n-1$ nonempty Young diagrams in a $1\times(n-1)$ grid,
and they are forced to be all rectangles. All this suggests that the rectangular torus chart $T_{\cC^R}$ should be the analogue of
the chart of the Landau-Ginzurg mirror corresponding to local systems on the Lagrangian Clifford torus in the projective
space, and we make this analogy precise in Theorem \ref{Thm1}.
\par
We now define more carefully the rectangular torus chart $T_{\cC^R}$. Recall from Section \ref{Section1}
the dual Pl\"{u}cker embedding
$$\Gr(n-k,n)\to\mathbb{P}^{{n \choose n-k}-1}\quad , \quad [\check{M}] \mapsto [\check{p}_d(\check{M})]$$
with the Pl\"{u}cker coordinates $\check{p}_{d}$ on the target parametrized by
$d\subseteq k\times (n-k)$, ordered by lexicographic order on their sets of
horizontal steps $d^-\subset [n]$. For $1\leq i\leq k$ and $1\leq j\leq n-k$ denote
$d_{i\times j}$ the rectangular $i\times j$ Young diagram, and $z_{ij}$ standard coordinates
on $(\mathbb{C}^\star)^{k(n-k)}$ (as usual in lexicographic order on the subscript). Define a map
$$\iota_{\cC^R}:(\mathbb{C}^\star)^{k(n-k)}\to\vGr(k,n)\quad , \quad \check{p}_d=
\begin{cases}
1 &\text{if } d=\emptyset\\
z_{ij} &\text{if } d=d_{i\times j}
\end{cases}\quad .
$$
This is well defined because the missing Pl\"{u}cker coordinates in the definition, i.e. those
$\check{p}_d$ that don't correspond to $d$ rectangular Young diagram, are determined by the
Pl\"{u}cker relations in $\mathbb{C}[\vGr(k,n)]$. Indeed, one can show that the rectangular
coordinates form a transcendence basis of the function
field $\mathbb{C}(\vGr(k,n))$, so that the non-rectangular Pl\"{u}cker coordinates are
rational functions of the rectangular ones; moreover, these rational functions must be Laurent polynomials.
This is an instance of the \textit{Laurent phenomenon} in the theory of cluster algebras, see for example \cite{RW}.
\par
The map above is an open embedding with image
$$T_{\cC^R}=\iota_{\cC^R}((\mathbb{C}^\star)^{k(n-k)})=\{\, [\check{M}]\in\vGr(k,n)\, :\, \check{p}_d([\check{M}])\neq 0\, \forall d\,\, \text{rectangular}\,\}$$
and this is what we call rectangular chart of $\vGr(k,n)$. The restriction
of the global function $W$ to the rectangular chart $T_{\cC^R}$ has been computed explicitly
by Marsh-Rietsch. We record their formula here for later use in Theorem \ref{Thm1}.

\begin{theorem}\label{ThmRM}(\cite[Section 6.3]{MR}; see also \cite[Proposition 9.5]{RW})
If $W:\vGr(k,n)\to\mathbb{C}$ is the Landau-Ginzburg superpotential, using the notation
$\check{p}_{i\times j}$ for $\check{p}_{d_{i\times j}}$ its restriction to the rectangular chart $T_{\cC^R}$ is given by
$$W\restriction_{T_{\cC^R}}=\check{p}_{1\times 1} +
\sum_{i=2}^{k}\sum_{j=1}^{n-k}\frac{   \check{p}_{i\times j}\check{p}_{(i-2)\times (j-1)}   }{   \check{p}_{(i-1)\times (j-1)}\check{p}_{(i-1)\times j}   }
+ \frac{   \check{p}_{(k-1)\times (n-k-1)}   }{   \check{p}_{k\times (n-k)}   } +
\sum_{i=1}^{k}\sum_{j=2}^{n-k}\frac{ \check{p}_{i\times j}\check{p}_{(i-1)\times (j-2)} }{ \check{p}_{(i-1)\times (j-1)}\check{p}_{i\times (j-1)} }
\quad .$$
\end{theorem}

We conclude this section with a criterion to decide when a critical
point of $W$ belongs to the rectangular chart $T_{\cC^R}$, formulated in terms of zeros of Schur
polynomials at roots of $\pm 1$.

\begin{proposition}\label{Prop2}
Let $I^\vee$ be a size $n-k$ subset of roots of $x^n=(-1)^{n-k+1}$ and denote $[\check{M}_{I^\vee}]\in\vGr(k,n)\subset\Gr(n-k,n)$
the corresponding critical point of $W:\vGr(k,n)\to\mathbb{C}$; then the following properties hold:
\begin{enumerate}
	\item $$[\check{M}_{I^\vee}]\in T_{\cC^R} \iff S_{d^T}(I^\vee)\neq 0 \quad \forall d\,\,\text{rectangular} $$
		where $S_{d^T}(x_1,\ldots ,x_{n-k})$ is the Schur polynomial of the transpose diagram of $d$,
		thought as diagram in the $(n-k)\times k$ grid.
	\item The dihedral group $D_n$ acts on the sets $I^\vee$ via $rI^\vee=e^{2\pi i/n}I^\vee$ and $sI^\vee=\overline{I^\vee}$
	and
	$$S_{d^T}(I^\vee)\neq 0 \iff S_{d^T}(gI^\vee)\neq 0\quad\forall g\in D_n\quad .$$
\end{enumerate}

\end{proposition}

\begin{proof}
\textit{1)} From the definition of $\iota_{\cC^R}$ with image the rectangular chart $T_{\cC^R}$
$$[\check{M}_{I^\vee}]\in T_{\cC^R} \iff \check{p}_d([\check{M}_{I^\vee}])\neq 0\quad \forall d \quad \text{rectangular Young diagram} \quad ,$$
and from the description of the critical points of $W$ (Theorem \ref{ThmKa}) given at the beginning
of this section
$$[\check{M}_{I^\vee}]=
\begin{bmatrix}
    1 & 1 & 1 & \dots  & 1 \\
    \zeta_1 & \zeta_2 & \zeta_3 & \dots  & \zeta_{n-k} \\
    \zeta_1^2 & \zeta_2^2 & \zeta_3^2 & \dots  & \zeta_{n-k}^2 \\
    \vdots & \vdots & \vdots &  & \vdots \\
    \zeta_1^{n-1} & \zeta_2^{n-1} & \zeta_3^{n-1} & \dots  & \zeta_{n-k}^{n-1}
\end{bmatrix}$$
where $I^\vee=\{\zeta_1,\ldots ,\zeta_{n-k}\}$ is a set of $n-k$ distinct roots of $x^n=(-1)^{n-k+1}$.
The horizontal steps of the empty diagram $d=\emptyset$ are $d^-=\emptyset^-=\{1,\ldots ,n-k\}$, so that
by the Vandermonde formula and the definition of $\check{p}_d$ as determinant of minor at rows $d^-$

$$\check{p}_\emptyset([\check{M}_{I^\vee}])=\prod_{1\leq i<j\leq n-k}(\zeta_j-\zeta_i)\neq 0 \quad .$$

Therefore

$$[\check{M}_{I^\vee}]\in T_{\cC^R} \iff \frac{\check{p}_d([\check{M}_{I^\vee}])}{\check{p}_{\emptyset}([\check{M}_{I^\vee}])}\neq 0\quad \forall d \quad \text{rectangular Young diagram}\quad .$$

The claim follows from the fact that

$$\frac{\check{p}_d([\check{M}_{I^\vee}])}{\check{p}_\emptyset([\check{M}_{I^\vee}])}=S_{d^T}(I^\vee)\quad .$$

This holds because when we write the diagram $d^T$ in partition form $d^T=(d^T_1,\ldots ,d^T_{n-k})$, i.e.
as a tuple where $d^T_i$ is the number of boxes at row $i$, one of the many equivalent ways
of defining the Schur polynomial of $d^T$ is
$$S_{d^T}(x_1,\ldots ,x_{n-k})=\frac{\operatorname{det}
	\begin{pmatrix}
    x_1^{d^T_{n-k}} & \dots & x_{n-k}^{d^T_{n-k}} \\
    x_1^{d^T_{n-k-1}+1} & \dots & x_{n-k}^{d^T_{n-k-1}+1} \\
    \vdots &  & \vdots \\
    x_1^{d^T_1+n-k-1} & \dots & x_{n-k}^{d^T_1+n-k-1}
	\end{pmatrix}}{\operatorname{det}
	\begin{pmatrix}
	1 & \dots & 1 \\
    x_1 & \dots & x_{n-k} \\
    \vdots &  & \vdots \\
    x_1^{n-k-1} & \dots & x_{n-k}^{n-k-1}
	\end{pmatrix}}\quad .$$
Once evaluated at $I^\vee$, the denominator is $\check{p}_\emptyset([\check{M}_{I^\vee}])$ and
the numerator equals $\check{p}_d([\check{M}_{I^\vee}])$ because the horizontal steps of $d$ and
the number of boxes in each row of $d^T$ are related by
$$d^-=\{d^T_{n-k}+1,\ldots ,d^T_1+n-k\}\quad .$$
\textit{2)} Observe that for every $d$
$$S_{d^T}(rI^\vee)=(e^{2\pi i/n})^{|{d^T}|}S_d(I^\vee) \quad ,\quad S_{d^T}(sI^\vee)=\overline{S_{d^T}(I^\vee)}$$
because $S_{d^T}$ is a homogeneous polynomial of degree the number of boxes $|d^T|$ of $d^T$.
\end{proof}

\section{Gelfand-Cetlin torus}\label{Section3}

From the point of view of symplectic topology, $\Gr(k,n)$ fits naturally in the class of
$\textit{coadjoint orbits}$. We denote $\mathfrak{u}(n)^\vee$ the dual Lie algebra of the
unitary group $U(n)$
and recall that the Lie bracket induces a symplectic structure on coadjoint orbits given by
$$\omega_\phi(X,Y) = \phi([X,Y])\quad \quad \phi \in \mathfrak{u}(n)^\vee\quad X,Y\in\mathfrak{u}(n) \quad ,$$
and the action of $U(n)$ is Hamiltonian with respect to this
structure. It is convenient to identify $\mathfrak{u}(n)^\vee$ with the real
vector space $\mathcal{H}_n$ of Hermitian matrices of size $n$, where the action of $U(n)$ is
given by conjugation. Each $H\in\mathcal{H}_n$ has real eigenvalues, labeled $\alpha_1\geq \ldots \geq\alpha_n$ .
By the spectral theorem, the tuple $\alpha = (\alpha_1, \ldots ,\alpha_n)$ is a complete
invariant for the orbits of the $U(n)$ action, so that we can denote them $\mathcal{O}_\alpha\subset\mathcal{H}_n$.
We are interested in the case where
$$\alpha_1=\ldots =\alpha_k > \alpha_{k+1} = \ldots = \alpha_n$$
because in this case $\mathcal{O}_\alpha \cong \Gr(k,n)$. Therefore, in this
section we think $\Gr(k,n)\subset \mathcal{H}_n$ as Hermitian matrices with $k$ equal big eigenvalues
and $n-k$ equal small eigenvalues.
\par
\begin{figure}[H]
  \centering
   %add desired spacing between images, e. g. ~, \quad, \qquad, \hfill etc. 
    %(or a blank line to force the subfigure onto a new line)
	\includegraphics[width=0.4\textwidth]{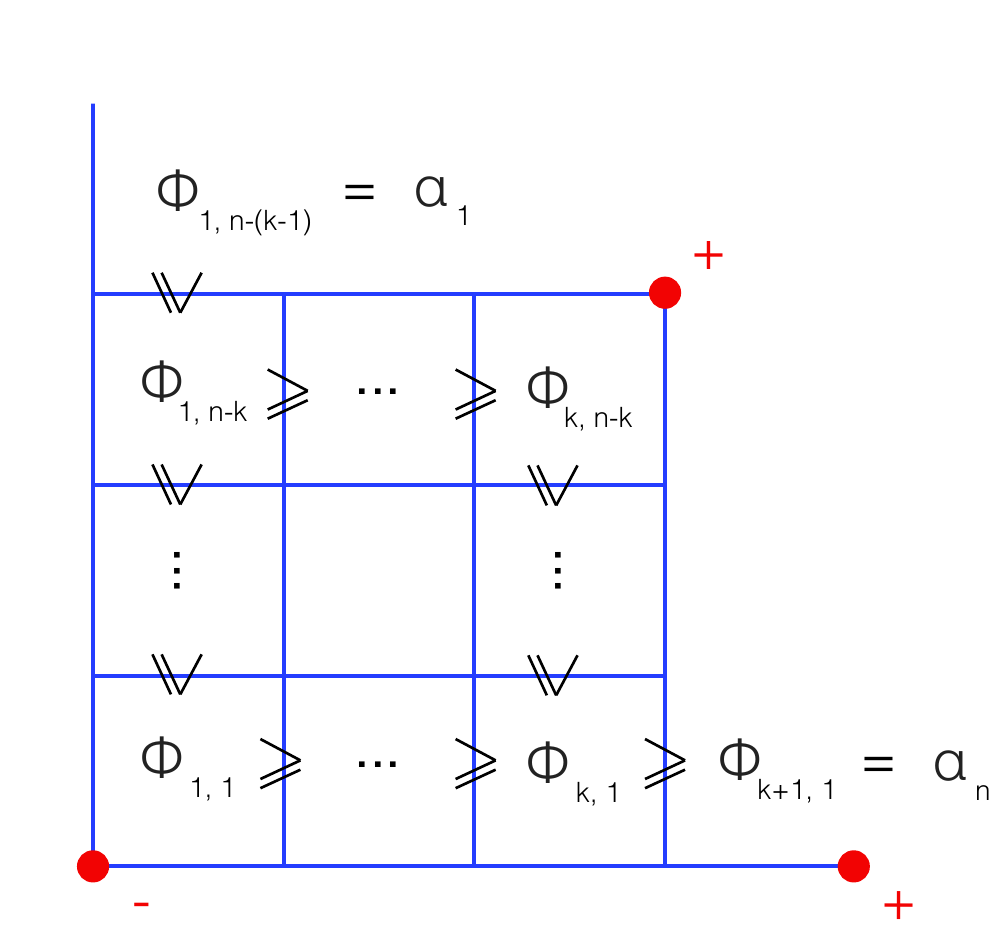}
    \caption{Inequalities for the Gelfand-Cetlin polytope $\Delta_\alpha$ and ladder diagram $\Gamma_\alpha$.}
    \label{fig:GC}
\end{figure}

For each $H\in\Gr(k,n)$ and $1\leq s \leq n-1$ one can consider the size $s$ minor $H_s$ of $H$ consisting
of the first $s$ rows and columns, so that $H_s\in\mathcal{H}_s$ and in analogy with what was done
above one can label its eigenvalues
$$\Phi_{1,s}\geq \Phi_{2,s-1}\geq \ldots\geq \Phi_{s-1,2}\geq \Phi_{s,1}\quad .$$
As a consequence of the min-max characterization of the eigenvalues,
for each $s\geq 2$ the eigenvalues of $H_s$ interlace with those of $H_{s-1}$:
\par
$$
\begin{matrix}
\Phi_{1,s}	&			& 				&			&\Phi_{2,s-1}\quad \dots	& \Phi_{s-1,2} & &				&			& \Phi_{s,1} \\
			& \brgeq 	& 				& \trgeq	&				&			 & \brgeq &				& \trgeq	&	\\
			&			& \Phi_{1,s-1}	&			&\				& \dots		 & & \Phi_{s-1,1}	&			&
\end{matrix}\quad .$$\\\\
This implies, together with the assumption that $\alpha_1=\ldots =\alpha_k$ and $\alpha_{k+1} = \ldots = \alpha_n$,
the inequalities of Figure \ref{fig:GC}.

\iffalse
$$
\begin{matrix}
    \Phi_{1,n-(k-1)}=\alpha_1 &  &  &  & \\
     &  &  &  &  \\
    \vgeq &  &  &  &  \\
     &  &  &  &  \\
    \Phi_{1,n-k} & \geq & \dots & \geq  & \Phi_{k,n-k} \\
     &  &  &  &  \\
    \vgeq &  &  &  & \vgeq \\
     &  &  &  &  \\
    \vdots &  &  &  & \vdots \\
     &  &  &  &  \\
    \vgeq &  & & & \vgeq \\
     &  &  &  &  \\
    \Phi_{1,1} & \geq & \dots & \geq  & \Phi_{k,1} & \geq & \Phi_{k+1,1}=\alpha_n

\end{matrix}$$\\\\
\fi
Therefore we are left with $k(n-k)$ nonconstant functions on $\Gr(k,n)$, and ordering
them by lexicographic order on the subscripts they become the entries of a map
$$\Phi : \Gr(k,n) \to \mathbb{R}^{k(n-k)}\quad .$$

The image $\Delta_\alpha$ of this map is a convex polytope
cut out by the inequalities of Figure \ref{fig:GC}. In fact, the following holds.
\begin{theorem}(Guillemin-Sternberg \cite[Section 5]{GS})
The map $\Phi$ is a completely integrable system on $\Gr(k,n)$: on the open dense subset
$\Phi^{-1}(\Delta_\alpha \setminus \partial\Delta_\alpha)\subset\Gr(k,n)$ its entries are smooth
functions, they pairwise Poisson commute, and their differentials are linearly independent.
\end{theorem}

It follows from general properties of completely integrable systems that the fibers of $\Phi$ over the interior of $\Delta_\alpha$ must be Lagrangian tori (see
for example Duistermaat \cite[Theorem 1.1]{Du}). Moreover, the following
result allows to locate the unique monotone fiber among these tori.

\begin{theorem}(Cho-Kim \cite[Section 5 and Theorem 5.2]{CK})
If $\alpha_1=\cdots = \alpha_k = n-k$ and $\alpha_{k+1}=\cdots = \alpha_n = -k$,
the symplectic structure $\omega$ on $\mathcal{O}_\alpha\cong \Gr(k,n)$ satisfies $[\omega]=c_1$
and the interior of $\Delta_\alpha$ contains a unique lattice point $\mathbf{u}$
such that $\Phi^{-1}(\mathbf{u})$ is monotone. The coordinates of this lattice point
are $\mathbf{u}_{i,j} = j-i$.
\end{theorem}
More generally, Cho-Kim \cite{CK} classify all the monotone Lagrangian
fibers of the Gelfand-Cetlin system for arbitrary partial flag manifolds.
We will denote the monotone torus fiber for the Grassmannian by $T^{k(n-k)}\subset\Gr(k,n)$, and
simply call it the Gelfand-Cetlin torus from now on.
\par
Being a regular Lagrangian fibration over the interior of $\Delta_\alpha$, the integrable system
$\Phi$ induces a basis $\gamma_{i,j}\in H_1(T^{k(n-k)};\mathbb{Z})$. Denote
$$W_{T^{k(n-k)}}\in\mathbb{C}[x_{i,j}^{\pm 1}] \quad\quad 1\leq i\leq k, \quad 1\leq j \leq n-k$$
the Maslov 2 disk potential written in the coordinates induced by this basis, as explained in
the Setup part. Our goal now is to compute this Laurent polynomial explicitly.
\par
The two key ingredients of this calculation are a correspondence between Maslov 2 $J$-holomorphic
disks and codimension 1 faces of $\Delta_\alpha$ due to Nishinou-Nohara-Ueda \cite{NNU}, and a combinatorial
description of boundary faces of $\Delta_\alpha$ due to An-Cho-Kim \cite{ACK}.

\begin{theorem}\label{ThmNNU}(Nishinou-Nohara-Ueda \cite[Theorem 10.1]{NNU})
The Maslov 2 disk potential $W_{T^{k(n-k)}}$ has one monomial with coefficient one for each
codimension 1 face of $\Delta_\alpha$, with exponents the coordinates of the corresponding
inward primitive normal vector.
\end{theorem}

The idea of the theorem above is to consider a toric degeneration of $\Gr(k,n)$ to a singular
toric variety $X(\Delta_\alpha)$ whose polytope is $\Delta_\alpha$, and to construct a \textit{small}
toric resolution of singularities for it. Small here means that the exceptional locus has codimension at
least two. Such degenerations have been found by Gonciuela-Lakshmibai \cite{GL}.
One then uses the degeneration to construct a cobordism between the moduli space of Maslov 2 disks bounding $T^{k(n-k)}\subset\Gr(k,n)$
and the moduli space of Maslov 2 disks bounding a toric Lagrangian fiber of the resolution of $X(\Delta_\alpha)$.
Smallness of the resolution guarantees that disks intersecting the singular locus in the
central fiber of the degeneration don't contribute to the potential, because they correspond to disks of
Maslov index at least 4 in the resolution. One then concludes by using the calculation of
Cho \cite{Ch} of the disk potential of toric Lagrangian fibers.
\par
In the proof of Proposition \ref{Prop3} we also make use of the following result.

\begin{theorem}\label{ThmACK}(An-Cho-Kim \cite[Theorem 1.11]{ACK})
There is a bijection between faces $\Delta_f\subset\Delta_\alpha$ and face graphs $\Gamma_f\subset\Gamma_\alpha$
in the ladder diagram, such that the dimension of $\Delta_f$ matches the number of
loops in $\Gamma_f$.
\end{theorem}

The \textit{ladder diagram} $\Gamma_\alpha$ is a $k\times (n-k)$ grid with
two extra edges; the bottom left corner is labelled by $-$ and there are two nodes labelled by $+$,
see Figure \ref{fig:GC}. A \textit{positive path} is
a sequence of edges in $\Gamma_\alpha$ that connects the $-$ node with one of the two $+$ nodes
and only goes up or right. A \textit{face graph} $\Gamma_f\subset\Gamma_\alpha$ is any union of positive
paths that covers both $+$ nodes. A loop of  $\Gamma_f$ is a minimal unoriented
cycle of $\Gamma_f$ thought as a graph.
\par
Figure \ref{fig:faces} gives a complete list of face graphs with five loops
for $k=2$ and $n=5$, and the bijection mentioned in the theorem above is obtained by setting
to $=$ those inequalities defining the polytope $\Delta_\alpha$ that do not cross an edge when
we put $\Gamma_f$ on top of the grid of inequalities as done in Figure \ref{fig:GC}.
\par
We use the above mentioned toric degeneration result of Nishinou-Nohara-Ueda
\cite{NNU} to work out an explicit formula in terms of $k$ and $n$ for the
disk potential of the Gelfand-Cetlin torus. This formula will be used in the
proof of Theorem \ref{Thm1}.

\begin{proposition}\label{Prop3}
If $T^{k(n-k)}\subset\Gr(k,n)$ is the Gelfand-Cetlin torus:
\begin{enumerate}
	\item There are $(k-1)(n-k)+k(n-k-1)+2$ Maslov 2 $J$-holomorphic disks through a generic
		point of the torus.
	\item The Maslov 2 disk potential is given by
	$$ W_{T^{k(n-k)}} = \sum_{i=1}^{k-1}\sum_{j=1}^{n-k}\frac{x_{i,j}}{x_{i+1,j}} +
	\sum_{i=1}^{k}\sum_{j=1}^{n-k-1}\frac{x_{i,j+1}}{x_{i,j}} + \frac{1}{x_{1,n-k}} + x_{k,1}\quad .$$
\end{enumerate}
\end{proposition}

\ytableausetup{boxsize=0.5em}

\begin{proof}
\textit{1)} Combining Theorems \ref{ThmNNU} and \ref{ThmACK} we know that through a
generic point of $T^{k(n-k)}$ there must be exactly one Maslov 2 $J$-holomorphic disk for
each codimension one face $\Delta_f\subset \Delta_\alpha$, and these in turn correspond to face graphs
$\Gamma_f\subset\Gamma_\alpha$ with $k(n-k)-1$ loops, where $k(n-k)=\operatorname{dim}\Delta_\alpha$.
The ambient ladder diagram $\Gamma_\alpha$ is a grid with $k(n-k)$ cells, therefore $\Delta_f$
falls in one of the following three types (see Figure \ref{fig:faces}):
\begin{itemize}
	\item \textbf{Type \ydiagram{2,0}} - the full $\Gamma_\alpha$ minus an interior vertical edge ;
	\item \textbf{Type \ydiagram{1,1}} - the full $\Gamma_\alpha$ minus an interior horizontal edge ;
	\item \textbf{Type \ydiagram{1,0}} - one of two exceptional cases consisting of the full $\Gamma_\alpha$ minus
		a loop at a corner .
\end{itemize}

\begin{figure}[H]
	\hspace{0.5cm}
   %add desired spacing between images, e. g. ~, \quad, \qquad, \hfill etc. 
    %(or a blank line to force the subfigure onto a new line)
    \begin{subfigure}[b]{0.15\textwidth}
        \includegraphics[width=\textwidth]{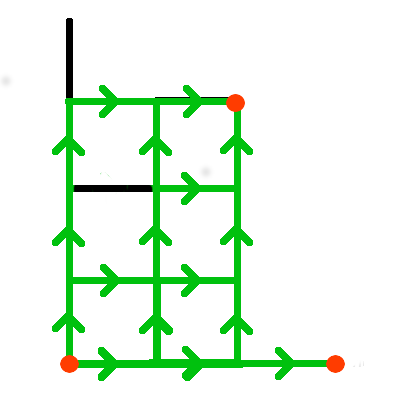}~
        \includegraphics[width=\textwidth]{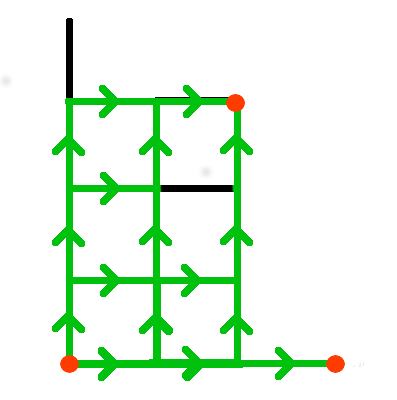}\quad
        \includegraphics[width=\textwidth]{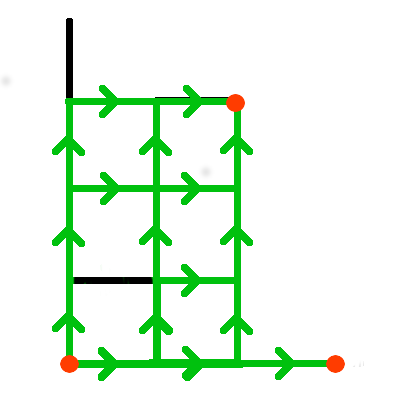}~
        \includegraphics[width=\textwidth]{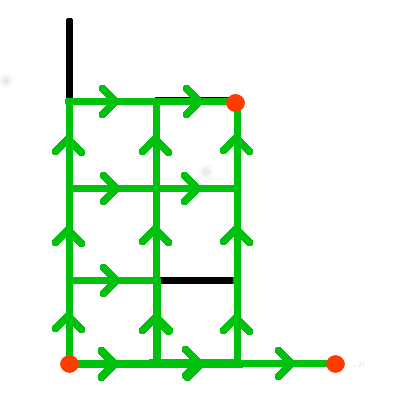}
        \caption{Type \ydiagram{1,1}}
    \end{subfigure}
    \hspace{3cm}
    %add desired spacing between images, e. g. ~, \quad, \qquad, \hfill etc. 
    %(or a blank line to force the subfigure onto a new line)
    \begin{subfigure}[b]{0.15\textwidth}
        \includegraphics[width=\textwidth]{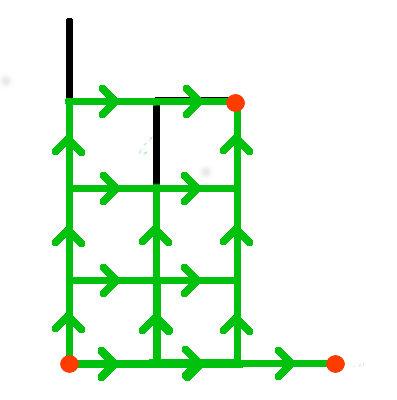}~
        \includegraphics[width=\textwidth]{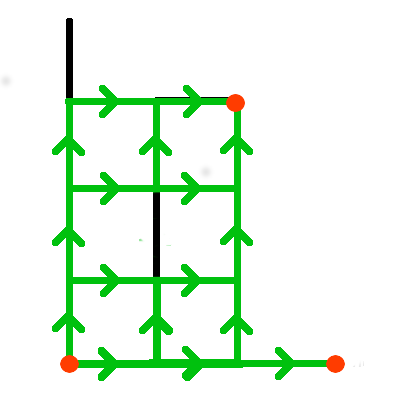}\quad
        \includegraphics[width=\textwidth]{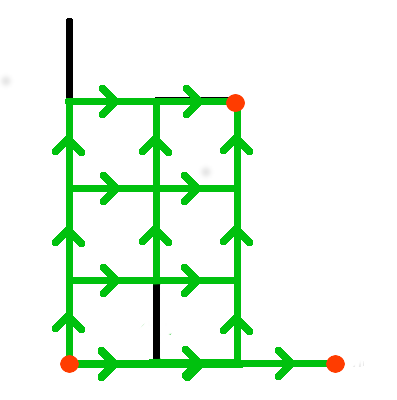}
        \caption{Type \ydiagram{2,0}}
    \end{subfigure}
    \hspace{3cm}
    %add desired spacing between images, e. g. ~, \quad, \qquad, \hfill etc. 
    %(or a blank line to force the subfigure onto a new line)
    \begin{subfigure}[b]{0.15\textwidth}
        \includegraphics[width=\textwidth]{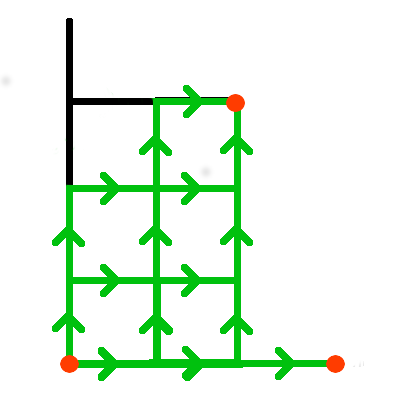}\quad
        \includegraphics[width=\textwidth]{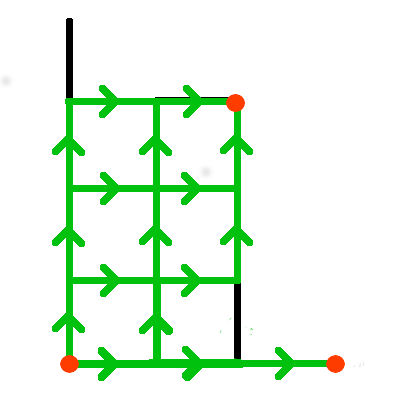}
        \caption{Type \ydiagram{1,0}}
    \end{subfigure}

    \caption{Face diagrams $\Gamma_f\subset\Gamma_\alpha$ of codimension 1 faces for the Gelfand-Cetlin
    	polytope of $\Delta_\alpha$ of the Grassmannian $\Gr(2,5)$.}
    \label{fig:faces}
\end{figure}

There are $(k-1)(n-k)$ ways to remove an interior vertical edge from the $k\times (n-k)$ grid $\Gamma_\alpha$:
they correspond to the ways of placing a horizontal brick \ydiagram{2,0}\, in it, where we think at
the central vertical edge as the one to be removed from $\Gamma_\alpha$ to obtain $\Gamma_f$. Similarly,
there are $k(n-k-1)$ ways to remove an interior horizontal edge from the $k\times (n-k)$ grid $\Gamma_\alpha$:
they correspond to the ways of placing a vertical brick \ydiagram{1,1}\, in it, where we think at
the central horizontal edge as the one to be removed from $\Gamma_\alpha$ to obtain $\Gamma_f$. Finally,
the two exceptional cases are best described by Figure \ref{fig:faces}. Note that one cannot get any
face graph $\Gamma_f$ by removing boundary edges at the top right and bottom left corners of
$\Gamma_\alpha$ without violating the condition on $\Gamma_f$ of being union of positive paths
and covering both the $+$ nodes.
\par
\textit{2)} We write down inequalities for the
codimension one faces $\Delta_f\subset \Delta_\alpha$ corresponding to the face graphs $\Gamma_f$ above,
and work out the coordinates of the corresponding inward primitive normal vectors $n_f\in\mathbb{Z}^{k(n-k)}\subset\mathbb{R}^{k(n-k)}$;
these coordinates will be the exponents of the Laurent monomials of $W_{T^{k(n-k)}}$, where the variables $x_{i,j}$ correspond
to the coordinates $\Phi_{i,j}$ of the Gelfand-Cetlin integrable system as explained at the beginning
of this section.
\par
By putting different types of face graphs $\Gamma_f$ on top of the grid of inequalities
of $\Delta_\alpha$ as in Figure \ref{fig:GC}, one finds:

\begin{itemize}
	\item \textbf{Type \ydiagram{2,0}} $\Phi_{i,j} = \Phi_{i+1,j}$ for $1\leq i\leq k-1$ and $1\leq j\leq n-k$ ;
	\item \textbf{Type \ydiagram{1,1}} $\Phi_{i,j} = \Phi_{i,j+1}$ for $1\leq i\leq k$ and $1\leq j\leq n-k-1$ ;
	\item \textbf{Type \ydiagram{1,0}} $\Phi_{1,n-k} = \alpha_1$ and $\Phi_{k,1} = \alpha_n$ .
\end{itemize}

Any of the codimension 1 faces $\Delta_f\subset \Delta_\alpha$ above is obtained by setting
one $\geq$ to $=$ in the facet presentation of the Gelfand-Cetlin polytope:
$$\Delta_\alpha =\{ \; \Phi\in\mathbb{R}^{k(n-k)} \; : \; n_f\cdot \Phi \geq -c_f \; \forall \Delta_f\subset\Delta_\alpha \; \textrm{such that} \; \operatorname{codim}\Delta_f=1 \; \}\quad .$$
Here $\cdot$ denotes the inner product
with the unique vector $n_f\in\mathbb{Z}^{k(n-k)}\subset\mathbb{R}^{k(n-k)}$
that generates the semigroup of integral vectors satisfying the inequality. In our case the inward primitive normal vectors are:

\begin{itemize}
	\item \textbf{Type \ydiagram{2,0}} for $1\leq i\leq k-1$ and $1\leq j\leq n-k$
		$$n_f = \begin{cases}
				1 &\text{in entry } \Phi_{i,j}\\
				-1 &\text{in entry } \Phi_{i+1,j}\\
				0 &\text{ in other entries}
				\end{cases} \quad ;$$
	\item \textbf{Type \ydiagram{1,1}} for $1\leq i\leq k$ and $1\leq j\leq n-k-1$
		$$n_f = \begin{cases}
				-1 &\text{in entry } \Phi_{i,j}\\
				1 &\text{in entry } \Phi_{i,j+1}\\
				0 &\text{ in other entries}
				\end{cases} \quad ;$$
	\item \textbf{Type \ydiagram{1,0}} for $\Phi_{1,n-k} = \alpha_1$ and $\Phi_{k,1} = \alpha_n$ respectively
		\begin{center}
		$n_f = \begin{cases}
				-1 &\text{in entry } \Phi_{1,n-k}\\
				0 &\text{ in other entries}
				\end{cases}$
							\quad	and \quad
		$n_f = \begin{cases}
				1 &\text{in entry } \Phi_{k,1}\\
				0 &\text{ in other entries}
				\end{cases}$ \quad .
		\end{center}
\end{itemize}

We conclude that $W_{T^{k(n-k)}}$ is a sum of monomials whose exponents are given by the coordinates
of the inward primitive normal vectors, giving the formula of the statement.

\end{proof}

\section{Main theorems}\label{Section4}

In this section we give proofs for the main theorems of this article.
The key diagram to have in mind is the following:

\begin{equation}\label{diag:1}
\begin{tikzcd}[row sep=2cm, column sep=2cm]
^{D_n \circlearrowright}(\mathbb{C}^\times)^{k(n-k)} \arrow[rd, "W_{T^{k(n-k)}}" left] \arrow[r, hook, "\theta_R"] & \vGr(k,n)^{\circlearrowleft D_n}\arrow[d, "W"] \\
	& \mathbb{C}^{\circlearrowleft D_n}
\end{tikzcd}
\end{equation}

Recall that $W$ is the function on the Landau-Ginzburg mirror described in Section \ref{Section1}:
\ytableausetup{boxsize=0.3em}
$$W = \frac{\check{p}_1^{\ydiagram{1}}}{\check{p}_1} + \ldots + \frac{\check{p}_n^{\ydiagram{1}}}{\check{p}_n} \quad .$$
$W_{T^{k(n-k)}}$ is the Maslov 2 disk potential of the Gelfand-Cetlin torus,
thought as an algebraic function on the space of local systems $(\mathbb{C}^\times)^{k(n-k)}$ and
computed in Section \ref{Section3}:
$$ W_{T^{k(n-k)}} = \sum_{i=1}^{k-1}\sum_{j=1}^{n-k}\frac{x_{i,j}}{x_{i+1,j}} +
	\sum_{i=1}^{k}\sum_{j=1}^{n-k-1}\frac{x_{i,j+1}}{x_{i,j}} + \frac{1}{x_{1,n-k}} + x_{k,1}\quad .$$
The target of both maps is $\mathbb{C}$, and it carries an action of the dihedral group $D_n$
by $2\pi /n$ counter-clockwise rotation and conjugation, which encodes the symmetries of the eigenvalues of $c_1\star$
on $\QH(\Gr(k,n))$ observed in Figure \ref{fig:flowers}.
\par
We will use an open embedding identifying the local systems of the Gelfand-Cetlin torus
with the rectangular chart $T_{\cC^R}$ of
the mirror $\vGr(k,n)$ defined in Section \ref{Section2}. The embedding $\theta_R:(\mathbb{C}^\times)^{k(n-k)}\to\vGr(k,n)$ is defined
by the equations
$$x_{i,j} = \frac{ \check{p}_{(k+1-i)\times j} }{ \check{p}_{(k-i)\times (j-1)} } \quad 1\leq i\leq k\quad 1\leq j\leq n-k \quad ,$$
with $\check{p}_{i\times j}=\check{p}_{d_{i\times j}}$ denoting the Pl\"{u}cker coordinate
corresponding to the $i\times j$ rectangular Young diagram.
\par
The definition above is phrased to be efficient for the purposes of Theorem \ref{Thm1}. One can use the equations to write
the coordinates $\check{p}_{i\times j}$ in terms of the coordinates $x_{i,j}$ on the space of local systems,
proceeding by lexicographic order on $i,j$. Since the $k(n-k)$ functions $\check{p}_{i\times j}$
give a transcendence basis for the function field $\mathbb{C}(\vGr(k,n))$, the other entries
$\check{p}_d$ of this map for $d$ non-rectangular Young diagram are determined by the Pl\"{u}cker relations.
The image of this embedding is the rectangular chart $T_{\cC^R}$ of Section \ref{Section2}
$$T_{\cC^R}=\{[\check{M}]\in\vGr(k,n)\, :\, \check{p}_d([\check{M}])\neq 0\, \forall d \, \text{rectangular}\,\}\quad .$$
A new ingredient in the diagram above is the action of the dihedral group $D_n$ on
the mirror $\vGr(k,n)$ (the action on the space of local systems $(\mathbb{C}^\times)^{k(n-k)}$ will
be defined as pull-back along $\theta_R$, see proof of Theorem \ref{Thm2}).

\begin{definition}\label{DefDihedralRep}
Indexing the standard basis $v_d$ of $\mathbb{C}^{{n \choose {k}}}$ with Young diagrams $d\subseteq k\times (n-k)$,
call dihedral projective representation the group morphism
$D_n\to \operatorname{PGL}({n \choose {n-k}},\mathbb{C})$
given on the generators $r,s$ of $D_n$ by
$$r\cdot v_d = (e^{2\pi i/n})^{|d|}v_d \quad ,\quad s\cdot v_d = v_{\operatorname{PD}(d)}$$
where $|d|$ denotes the number of boxes in $d$ and $\operatorname{PD}(d)$ is the Poincar\'{e} dual
Young diagram of $d$.
\end{definition}

The \textit{Poincar\'{e} dual} of a Young diagram $d$ is obtained by taking the complement
of $d$ in the ambient $k\times (n-k)$ grid and rotating by $\pi$ to place it in the top left corner.

\begin{example}
Let $k=3$ and $n=7$. An example of Young diagram and its Poincar\'{e} dual
in the $3\times 4$ grid is given by

\ytableausetup{boxsize=normal}
$$d = \quad \ydiagram{4,3,0}\quad \textrm{and} \quad \operatorname{PD}(d) = \quad \ydiagram{4,1,0} \quad .$$

\end{example}

Observe that $r^n=s^2=1$ in the definition above because $e^{2\pi i/n}$ is an $n$-th root of
unity and $\operatorname{PD}(\operatorname{PD}(d))=d$. Moreover, the relation $rs=sr^{-1}$ holds
in $\operatorname{PGL}({n \choose {n-k}},\mathbb{C})$ because
$$rs\cdot v_d = (e^{2\pi i/n})^{|\operatorname{PD}(d)|}v_{\operatorname{PD}(d)}\quad ,\quad
sr^{-1}\cdot v_d = (e^{2\pi i/n})^{-|d|}v_{\operatorname{PD}(d)}$$
and $|\operatorname{PD}(d)|= k(n-k) - |d|$, so that $rs = e^{2\pi ik(n-k)/n}sr^{-1}$ with
the scaling factor independent of $d$. The dihedral projective representation
induces an algebraic action of $D_n$ on $\mathbb{P}^{{n \choose {k}}-1}$,
and the following holds.

\begin{lemma}\label{lemma:dihedral-action}
The image of the dual Pl\"{u}cker embedding $\Gr(n-k,n)\subset\mathbb{P}^{{n \choose {k}}-1}$ is
invariant under the action of $D_n$ induced by the dihedral projective representation.
\end{lemma}

\begin{proof}
Write the full rank $n\times (n-k)$ matrix $\check{M}$
as a list of rows
$$\check{M} = (m_1, \ldots, m_n) \quad m_i\in\mathbb{R}^{n-k} \; 1\leq i \leq n \quad .$$
If $d\subseteq k\times (n-k)$, and $d^-=\{i_1,\ldots ,i_{n-k}\}$
are its horizontal steps, the corresponding dual Pl\"{u}cker coordinate of
$[\check{M}]\in\Gr(n-k,n)$ is
$$\check{p}_d([\check{M}])=\operatorname{det}(m_{i_1},\ldots , m_{i_{n-k}}) \quad .$$
If $\zeta = e^{2\pi i/n}$, one can define a $D_n$-action on $\Gr(n-k,n)$ via
$$r\cdot [\check{M}] = [\zeta m_1,\zeta^2 m_2, \ldots , \zeta^n m_n] \quad \textrm{and} \quad s\cdot [\check{M}] = [m_n, m_{n-1},\ldots , m_1] \quad ,$$
and observe that
$$\check{p}_d(r\cdot [\check{M}]) = \zeta^{i_1+\cdots +i_{n-k}}\check{p}_d([\check{M}]) \quad , \quad
\check{p}_d(s\cdot [\check{M}]) = \check{p}_{PD(d)}([\check{M}]) \quad .$$
We claim that these agree with the homogeneous coordinates of $[\check{M}]$
under the projective action $D_n$ coming from the dihedral projective
representation. For the $s$ action there
is nothing to check. Regarding $r$, the claim follows from the existence of
a constant $C_{k,n}$ independent of $d$ such that $|d| = C_{k,n} + i_1 + \ldots + i_{n-k}$.
To see that this constant exists, write the Young diagram as a partition
$d=(d_1,\ldots d_k)$, with $d_i$ number of boxes in the $i$-th row of $d$.
If $d^|=\{j_1, \ldots j_{k}\}$ are the vertical steps of $d$, one has
$$d_1 = n-k-(j_1-1)\quad , \quad d_2=d_1-(j_2-j_1-1)\quad ,\quad \ldots\quad ,\quad d_k = d_{k-1} - (j_k-j_{k-1}-1) \quad .$$
It follows that
$$|d| = d_1+\cdots +d_k = k(n-k)+\frac{k(k+1)}{2}-(j_1+\cdots + j_k) \quad ;$$
moreover $[n] = d^- \sqcup d^|$ implies
$$i_1+\cdots + i_{n-k} = \frac{n(n+1)}{2} - (j_1+\cdots + j_k) \quad ,$$
from which one finds
$$C_{k,n} = k(n-k) + \frac{k(k+1)}{2} - \frac{n(n+1)}{2} \quad .$$
\end{proof}

Recall that the Landau-Ginzburg mirror of $\Gr(k,n)$ is
$$\vGr(k,n)=\Gr(n-k,n)\;\setminus\;\{ \; \check{p}_1\cdots\check{p}_n=0 \; \} \quad ,$$
where $\check{p}_1, \ldots , \check{p}_n$ are the Pl\"{u}cker coordinates of the $n$ boundary
rectangular Young diagrams, i.e. rectangular with at least one side covering a full side of the $k\times (n-k)$ grid
(and the $\emptyset$ diagram). Since this collection of diagrams is closed under Poincar\'{e} duality,
the divisor is $D_n$-invariant and one gets an action of $D_n$ on $\vGr(k,n)$.

\begin{definition}\label{DefYoungAction}
Call Young action the algebraic $D_n$-action on $\vGr(k,n)$ induced by
the dihedral projective representation of Definition \ref{DefDihedralRep}.
\end{definition}

For notational convenience, we introduce a $D_n$-action on the sets $I$ of $k$ distinct
roots of $x^n=(-1)^{k+1}$ that extends the one on $\mathbb{C}$ element-wise
$$r\cdot I = e^{2\pi i/n}I \quad ,\quad s\cdot I = \overline{I}\quad .$$
We also denote $I^c$ the set of $n-k$ roots of $x^n=(-1)^{k+1}$ that are not
in $I$, and observe that the formula $I^\vee = e^{\pi i}I^c$ gives a bijection
between size $k$ sets of roots of $x^n=(-1)^{k+1}$ and size $n-k$ sets of
roots of $x^n=(-1)^{n-k+1}$. In the theorems below, $I$ will parametrize objects
supported on the Gelfand-Cetlin torus $T^{k(n-k)}\subset\Gr(k,n)$, whereas
$I^\vee$ will parametrize critical points of the Landau-Ginzburg potential
$W:\vGr(k,n)\to\mathbb{C}$.
\par
The values of Schur polynomials at roots in $I$ and $I^\vee$
are related by the following.

\begin{lemma}\label{DualSchur}(Rietsch \cite[Lemma 4.4]{Ri})
If $d\subseteq k\times (n-k)$, and $d^T\subseteq(n-k)\times k$ is the transpose
diagram, then
$$S_{d^T}(I^\vee) = S_d(I) \quad .$$
\end{lemma}

We are now ready to prove the following.

\begin{theorem}\label{Thm1}
Choosing $I_0$ to be the set of $k$ roots of $x^n=(-1)^{k+1}$ closest to 1, the objects
obtained by giving $T^{k(n-k)}$ the different local systems
$$T^{k(n-k)}_{I_0} \; ,\; T^{k(n-k)}_{rI_0} \; ,\; T^{k(n-k)}_{r^2I_0} \;, \; \ldots  \; ,\; T^{k(n-k)}_{r^{n-1}I_0}$$
are defined and split-generate the $n$ summands $\Fuk_\lambda(\Gr(k,n))$ of the monotone
Fukaya category with maximum $|\lambda|$.
\end{theorem}

\begin{proof}
As a first step we prove that diagram (\ref{diag:1}) commutes,
or in other words $\theta_R^*W=W_{T^{k(n-k)}}$.
The image of $\theta_R$ is the rectangular torus chart $T_{\cC^R}\subset\vGr(k,n)$
and by Theorem \ref{ThmRM}
$$W\restriction_{T_{\cC^R}}=\check{p}_{1\times 1} +
\sum_{i=2}^{k}\sum_{j=1}^{n-k}\frac{   \check{p}_{i\times j}\check{p}_{(i-2)\times (j-1)}   }{   \check{p}_{(i-1)\times (j-1)}\check{p}_{(i-1)\times j}   }
+ \frac{   \check{p}_{(k-1)\times (n-k-1)}   }{   \check{p}_{k\times (n-k)}   } +
\sum_{i=1}^{k}\sum_{j=2}^{n-k}\frac{ \check{p}_{i\times j}\check{p}_{(i-1)\times (j-2)} }{ \check{p}_{(i-1)\times (j-1)}\check{p}_{i\times (j-1)} }
\quad .$$
From the definition of $\theta_R$
$$x_{i,j} = \frac{ \check{p}_{(k+1-i)\times j} }{ \check{p}_{(k-i)\times (j-1)} } \quad 1\leq i\leq k\quad 1\leq j\leq n-k$$
so that
$$\theta_R^*W = x_{k,1} + \sum_{i=2}^{k}\sum_{j=1}^{n-k}\frac{ x_{k+1-i,j} }{ x_{k+2-i,j} } + \frac{1}{x_{1,n-k}} +
\sum_{i=1}^{k}\sum_{j=2}^{n-k}\frac{ x_{k+1-i,j} }{ x_{k+1-i,j-1} }\quad .$$
The last expression matches the Maslov 2 disk potential $W_{T^{k(n-k)}}$
computed in Proposition \ref{Prop3} (re-index the first sum with $l=k+1-i$, and
the second sum with $l=k+1-i$ and $h=j-1$).
\par
Because of Theorem \ref{Fact2}, critical points $\operatorname{hol}_I\in(\mathbb{C}^\times)^{k(n-k)}$
of $W_{T^{k(n-k)}}$ correspond to holonomies of local systems on the Gelfand-Cetlin torus
such that
$$\HF(T^{k(n-k)}_I, T^{k(n-k)}_I)\neq 0 \quad ,$$
and by commutativity of diagram (\ref{diag:1})
these are the same as critical points of $W:\vGr(k,n)\to\mathbb{C}$
that are contained in the rectangular chart $T_{\cC^R}$. By Theorem \ref{ThmKa} and
Lemma \ref{DualSchur}, critical points
of $W$ are of the form $[\check{M}_{I^\vee}]\in\vGr(k,n)$ for $I$ set of $k$ distinct roots
of $x^n=(-1)^{k+1}$, and by Proposition \ref{Prop2} $[\check{M}_{I^\vee}]\in T_{\cC^R}$ if and only
if the following nonvanishing condition on Schur polynomials holds:
$$S_{d^T}(I^\vee) = S_d(I) \neq 0 \quad \forall d \quad \text{rectangular}$$
(the first equality is Lemma \ref{DualSchur}).
As explained in the Introduction, this allows us to define for any $[\check{M}_{I^\vee}]\in T_{\cC^R}$ a local system on
the Gelfand-Cetlin torus with
$$\operatorname{hol}_{I}(\gamma_{ij})=\frac{S_{(k+1-i)\times j}(I)}{S_{(k-i)\times (j-1)}(I)}\quad .$$
Observe that arguing as in Proposition \ref{Prop2}
$$\frac{\check{p}_{(k+1-i)\times j}([\check{M}_{I^\vee}])}{\check{p}_{(k-i)\times (j-1)}([\check{M}_{I^\vee}])}
=\frac{S_{((k+1-i)\times j)^T}(I^\vee)}{S_{((k-i)\times (j-1))^T}(I^\vee)}
=\frac{S_{(k+1-i)\times j}(I)}{S_{(k-i)\times (j-1)}(I)}\quad .$$
We conclude that $\theta_R(\operatorname{hol}_I)=[\check{M}_{I^\vee}]$, meaning that
$\HF(T^{k(n-k)}_I,T^{k(n-k)}_I)\neq 0$.
\par
To decide in which summand $\Fuk_\lambda(\Gr(k,n))$ of the Fukaya category the
object $T^{k(n-k)}_I$ lives, we compute
\ytableausetup{boxsize=0.3em}
$$\lambda = W_{T^{k(n-k)}}(\operatorname{hol}_I)=W([\check{M}_{I^\vee}])=\left( \frac{\check{p}_1^{\ydiagram{1}}}{\check{p}_1} + \ldots + \frac{\check{p}_n^{\ydiagram{1}}}{\check{p}_n}\right)([\check{M}_{I^\vee}])\quad .$$
Dividing numerator and denominator by $\check{p}_\emptyset$, one finds for every $1\leq t \leq n$
$$\frac{ {\check{p}_t}^{ \ydiagram{1} } }{ \check{p}_t }([\check{M}_{{I}^\vee}]) = \frac{ S_{(\ydiagram{1}\star d_t)^T}(I^\vee) }{ S_{(d_t)^T}(I^\vee) }
=\frac{ S_{\ydiagram{1}\star d_t}(I) }{ S_{d_t}(I) } \quad ,$$
where $d_t$ is the $t$-th boundary rectangular Young diagram
(for the definition of boundary and interior rectangular diagrams, see Section \ref{Section1}).
\par
Recall from Proposition \ref{Prop1} that we called $\sigma_{d_t}\in\QH(\Gr(k,n))$ the Schubert class of
the Young diagram $d_t$, and Lemma \ref{lemma:schur-basis} says that
the action of $\sigma_t\star$  on the Schur basis $\sigma_I$ is given by
$$\sigma_{d_t}\star\sigma_I = S_{d_t}(I)\sigma_I\quad .$$
Since $d_t^{\ydiagram{1}}=\ydiagram{1}\star d_t$ is a single Schubert class due to the special
rectangular shape of $d_t$, we have
$$\sigma_{d_t^{\ydiagram{1}}}\star\sigma_I = S_{\ydiagram{1}}(I)S_{d_t}(I)\sigma_I$$
and therefore $S_{d_t^{\ydiagram{1}}}(I)=S_{\ydiagram{1}}(I)S_{d_t}(I)$. We conclude that
$\lambda = nS_{\ydiagram{1}}(I)$.
\par
Choosing now $I=I_0$ set of $k$ roots closest to $1$, the set $I_0^\vee$ contains the
$n-k$ roots of $x^n=(-1)^{n-k+1}$ closest to $1$, and the corresponding
critical point $[\check{M}_{I_0^\vee}]$ lies in the totally positive part of the Grassmannian;
see \cite[Theorem 1.1]{Ka}. This means that all the Pl\"{u}cker coordinates of
$[\check{M}_{I_0^\vee}]$ are real, nonvanishing, and have the same sign. In particular,
$[\check{M}_{I_0^\vee}]$ is in the rectangular chart $T_{\cC^R}$ and this gives
a nonzero object

$$T^{k(n-k)}_{I_0}\,\,\text{in}\,\,\Fuk_\lambda(\Gr(k,n)) \quad ,\quad \lambda = W([\check{M}_{I_0^\vee}])=nS_{\ydiagram{1}}(I_0) \in\mathbb{R}^+\quad .$$

Moreover, using the inequalities of Schur polynomials at roots of unity obtained by Rietsch \cite[Proposition 11.1]{Ri}
$$|S_{\ydiagram{1}}(I)|\leq S_{\ydiagram{1}}(I_0) \quad \forall I$$
and this says that $T^{k(n-k)}_{I_0}$ lives in a summand of the Fukaya category labelled by
an eigenvalue of $c_1\star$ with maximum modulus. It is known \cite[Proposition 3.3 and Corollary 4.11]{CL} that there are
$n$ eigenvalues of $c_1\star$ with maximum modulus, all of multiplicity 1, and they are given
by rotations of multiples of $2\pi /n$ of $nS_{\ydiagram{1}}(I_0)\in\mathbb{R}^+$. One can get
nonzero objects in each of them by considering
$$T^{k(n-k)}_{I_0},\, T^{k(n-k)}_{rI_0},\, T^{k(n-k)}_{r^2I_0},\, \ldots ,\, T^{k(n-k)}_{r^{n-1}I_0}$$
and observing that for all $0\leq s <n$ the critical point $[\check{M}_{(r^sI_0)^\vee}]$ of $W$
still belongs to the rectangular chart $T_{\cC^R}$, thanks to Proposition \ref{Prop2} (see also
proof of Theorem \ref{Thm2}), with
$$W([\check{M}_{(r^sI_0)^\vee}])=nS_{\ydiagram{1}}(r^sI_0)=r^{s}nS_{\ydiagram{1}}(I_0) \quad ,$$
where in the last step we used that the Schur polynomial of $\ydiagram{1}$ is linear. The fact
that these summands are indexed by eigenvalues of multiplicity one guarantees that a single
nonzero object generates, thanks to Theorem \ref{Fact3}.
\end{proof}

\begin{remark}\label{RemarkThm1}
In fact, the theorem above shows that whenever a critical point $[\check{M}_{I^\vee}]$ of the
Landau-Ginzburg superpotential belongs to the rectangular torus chart $T_{\cC^R}\subset\vGr(k,n)$
the object $T^{k(n-k)}_I$ is defined and nonzero in the summand of the Fukaya category $\Fuk_\lambda(\Gr(k,n))$ with
$$\lambda = W([\check{M}_{I^\vee}])=nS_{\ydiagram{1}}(I)$$
even when $|\lambda|$ has not maximum modulus. The difference in this case is that the
summand $\QH_\lambda(\Gr(k,n))$ of quantum cohomology is not necessarily one-dimensional,
therefore one cannot conclude that this object generates $\D\Fuk_\lambda(\Gr(k,n))$.
\end{remark}

As Figure \ref{fig:branes} illustrates, the question of what summands with lower $|\lambda|$
contain nonzero objects supported on the Gelfand-Cetlin torus is related to
the arithmetic of $k$ and $n$. On the other hand, we show in Theorem \ref{Thm2} that the equivariance
of the Landau-Ginzburg mirror with respect to the action of $D_n$ introduced at the beginning
of this section forces a certain dichotomy, for which summands $\Fuk_\lambda(\Gr(k,n))$
in each level $|\lambda|$ contain either none or a full $D_n$-orbit of nonzero objects.

\begin{theorem}\label{Thm2}
If $T^{k(n-k)}_I$ is an object of $\Fuk_\lambda(\Gr(k,n))$, then it is nonzero and the objects
$$T^{k(n-k)}_{gI}\quad \text{in}\quad \Fuk_{g\lambda}(\Gr(k,n)) \quad \text{for all}\quad g\in D_n$$
are defined and nonzero as well.
\end{theorem}

\begin{proof}
If $[\check{M}_{I^\vee}]\in\vGr(k,n)$ is a critical point of the Landau-Ginzburg
superpotential $W$, the Young action of $D_n$ on $\vGr(k,n)$ sends it to another
critical point
$$g[\check{M}_{I^\vee}]=[\check{M}_{(gI)^\vee}] \quad\quad \forall g\in D_n\quad .$$
We verify this equality on generators $r,s$ of $D_n$. Calling $N ={n \choose n-k}-1$, the Pl\"{u}cker
coordinates of $[\check{M}_{I^\vee}]\in\vGr(k,n)\subset\mathbb{P}^N$ are
$$[\check{p}_{d_0}([\check{M}_{I^\vee}]):\check{p}_{d_1}([\check{M}_{I^\vee}]):\cdots :\check{p}_{d_N}([\check{M}_{I^\vee}])]\quad .$$
The Young diagrams $d_0,\ldots ,d_N$ are labelled according to lexicographic order on
their sets of horizontal steps. The action of $r$ gives
$$r\cdot [\check{M}_{I^\vee}]=[(e^{2\pi i/n})^{|d_0|}\check{p}_{d_0}([\check{M}_{I^\vee}]):(e^{2\pi i/n})^{|d_1|}\check{p}_{d_1}([\check{M}_{I^\vee}]):\cdots :(e^{2\pi i/n})^{|d_N|}\check{p}_{d_N}([\check{M}_{I^\vee}])]$$
Now observe that $d_0=\emptyset$ is the empty diagram, whose horizontal steps are $\{1,\ldots n-k\}$.
Therefore the number of boxes of $d_0$ is $|d_0|=0$, and thanks to Proposition \ref{Prop2}
we know that $\check{p}_{d_0}([\check{M}_{I^\vee}])\neq 0$. Scaling the homogeneous coordinates by
this factor we have
$$r\cdot [\check{M}_{I^\vee}]=[1:(e^{2\pi i/n})^{|d_1|}S_{d_1}(I):\cdots :(e^{2\pi i/n})^{|d_N|}S_{d_N}(I)]\quad .$$
The Schur polynomial $S_{d_i}$ is homogeneous of degree $|d_i|$, therefore
$$(e^{2\pi i/n})^{|d_1|}S_{d_i}(I)=S_{d_i}(e^{2\pi i/n}I)=S_{d_i}(rI)\quad\quad \forall 0\leq i\leq N\quad .$$
Scaling back the homogeneous coordinates by a factor $\check{p}_{d_0}([\check{M}_{(rI)^\vee}])$ we get
$$r\cdot [\check{M}_{I^\vee}]=[\check{p}_{d_0}([\check{M}_{(rI)^\vee}]):\check{p}_{d_1}([\check{M}_{(rI)^\vee}]):\cdots :\check{p}_{d_N}([\check{M}_{(rI)^\vee}])]=[\check{M}_{(rI)^\vee}]\quad .$$
The verification for the action of the generator $s$ is analogous:
$$s\cdot [\check{M}_{I^\vee}]=[\check{p}_{PD(d_0)}([\check{M}_{I^\vee}]):\check{p}_{PD(d_1)}([\check{M}_{I^\vee}]):\cdots :\check{p}_{PD(d_N)}([\check{M}_{I^\vee}])]\quad .$$
Scaling the homogeneous coordinates by $\check{p}_{d_0}([\check{M}_{I^\vee}])$ we have
$$s\cdot [\check{M}_{I^\vee}]=[S_{PD(d_0)}(I):S_{PD(d_1)}(I):\cdots :S_{PD(d_N)}(I)]\quad .$$
From the interpretation of Schur polynomials $S_{d_i}$ as characters of representations of $GL(k,\mathbb{C})$,
it follows that
$$S_{PD(d_i)}(I)=\overline{S_{d_i}(I)}c_I \quad\quad \forall 0\leq i\leq N\quad .$$
Here $c_I=S_{d_N}(I)\in\mathbb{C}^\times$ is a constant that depends on $I$, but not on $d_i$ (see
for example of \cite[Lemma 4.4]{Ri}). Therefore we have
$$s\cdot [\check{M}_{I^\vee}]=[\overline{S_{d_0}(I)}:\overline{S_{d_1}(I)}:\cdots :\overline{S_{d_N}(I)}]=[S_{d_0}(sI):S_{d_1}(sI):\cdots :S_{d_N}(sI)]=[\check{M}_{(sI)^\vee}]$$
where the last equality is obtained by scaling the homogeneous coordinates by $\check{p}_{d_0}([\check{M}_{(sI)^\vee}])$.
This concludes the proof that the critical locus of $W$ is $D_n$-invariant.
\par
We also observe that $W$ is $D_n$-equivariant when restricted to the singular locus. Indeed,
we saw in Theorem \ref{Thm1} that $W([\check{M}_I])=nS_{\ydiagram{1}}(I)$ and therefore
$$r\cdot W([\check{M}_{I^\vee}])=e^{2\pi i/n}nS_{\ydiagram{1}}(I)=nS_{\ydiagram{1}}(rI)=W([\check{M}_{(rI)^\vee}])=W(r\cdot [\check{M}_{I^\vee}])$$
$$s\cdot W([\check{M}_{I^\vee}])=\overline{nS_{\ydiagram{1}}(I)}=nS_{\ydiagram{1}}(sI)=W([\check{M}_{(sI)^\vee}])=W(s\cdot [\check{M}_{I^\vee}])\quad .$$
Also notice that $W:\vGr(k,n)\to\mathbb{C}$ can't be globally $D_n$-equivariant,
because $W(s\cdot -)$ is an algebraic function while $s\cdot W(-)$ is not. On the other
hand $W$ is globally $\mathbb{Z}/n\mathbb{Z}$-equivariant, where $\mathbb{Z}/n\mathbb{Z}$ is
the subgroup of $D_n$ generated by $r$. Indeed, by definition we have
\ytableausetup{boxsize=0.3em}
$$W = \frac{\check{p}_1^{\ydiagram{1}}}{\check{p}_1} + \ldots + \frac{\check{p}_n^{\ydiagram{1}}}{\check{p}_n}$$
and recall that $\check{p}_t=\check{p}_{d_t}$ Pl\"{u}cker coordinate of the $t$-th boundary rectangular
Young diagram, with $\check{p}_t^{\ydiagram{1}}=\check{p}_{\ydiagram{1}\star d_t}$. The fact
that $|\ydiagram{1}\star d_t|=|d_t|+1$ for every $1\leq t\leq n$ implies that
$$\frac{r\cdot \check{p}_t^{\ydiagram{1}}}{r\cdot \check{p}_t}=\frac{(e^{2\pi i/n})^{|\ydiagram{1}\star d_t|}\check{p}_n^{\ydiagram{1}}}{(e^{2\pi i/n})^{|d_t|}\check{p}_n}=e^{2\pi i/n}\frac{\check{p}_n^{\ydiagram{1}}}{\check{p}_n}$$
which means that $W$ is globally $\mathbb{Z}/n\mathbb{Z}$-equivariant.
\par
Observe now that the rectangular torus chart $T_{\cC^R}\subset\vGr(k,n)$ is $\mathbb{Z}/n\mathbb{Z}$-invariant
because
$$\check{p}_d([\check{M}])\neq 0 \iff r\cdot \check{p}_d([\check{M}])=(e^{2\pi i/n})^{|d|}\check{p}_d([\check{M}])\neq 0\quad .$$
On the other hand $T_{\cC^R}$ is not $D_n$-invariant, because the Poincar\'{e} dual of an
interior rectangular Young diagram is not necessarily rectangular (as opposed to the case
of boundary rectangular Young diagrams whose duals are rectangular, see Section \ref{Section1}).
\par
Thanks to the fact that the action of $s$ on the critical points matches the action by
conjugation, we have
$$[\check{M}_{I^\vee}]\in T_{\cC^R} \iff [\check{M}_{I^\vee}]\in V \quad ,\quad V=\bigcap_{g\in D_n}g\cdot T_{\cC^R}$$
and $V$ is an open $D_n$-invariant subscheme of $T_{\cC^R}$. Defining $U=\theta_R^{-1}(V)\subset (\mathbb{C}^\times)^{k(n-k)}$, this
is an open subscheme of the space of rank one $\mathbb{C}$-linear local systems on the Gelfand-Cetlin
torus $T^{k(n-k)}\subset \Gr(k,n)$. Since $\theta_R$ is an open embedding, we can use it
to pull-back the $D_n$ action on $U$, so that diagram (\ref{diag:1}) becomes fully $\mathbb{Z}/n\mathbb{Z}$-equivariant,
and $D_n$ equivariant whenever restricted to critical points and values.
\end{proof}

In Theorem \ref{Thm3} we show that there are indeed classes of Grassmannians for which the considerations
of Theorem \ref{Thm1} and \ref{Thm2} suffice to identify a complete set of generators for the Fukaya category,
and prove homological mirror symmetry.

\begin{theorem}\label{Thm3}
When $n=p$ is prime the objects $T^{k(p-k)}_I$ split-generate the Fukaya category of $\Gr(k,p)$, and
for every $\lambda\in\mathbb{C}$ there is an equivalence of triangulated categories
$$\D\Fuk_\lambda(\Gr(k,p))\simeq \D\mathcal{S}(W^{-1}(\lambda))\quad .$$
\end{theorem}

\begin{proof}
By part (3) of Proposition \ref{Prop1} and the assumption $n=p$ prime, we have that
$$\operatorname{dim}\QH_\lambda(\Gr(k,p))=1 \quad\quad \forall \lambda \, \text{ eigenvalue of } c_1\star\quad .$$
Thanks to Theorem \ref{Fact3}, any nonzero object supported on the Gelfand-Cetlin torus will generate
the summand of the Fukaya category in which it lives, therefore it suffices to show that the torus supports objects in all summands to
have a complete set of generators.
\par
When $n=2$ we have $\Gr(1,2)=\mathbb{P}^1$, and the objects of the statement are the two local
system on the Clifford torus giving objects with nontrivial Floer cohomology investigated by
Cho \cite{Ch}. We will show that when $n=p>2$ is prime all critical points $[\check{M}_{I^\vee}]$ of $W:\vGr(k,p)\to\mathbb{C}$ are contained
in the rectangular chart $T_{\cC^R}$. By Theorem \ref{Thm1} and Remark \ref{RemarkThm1} this will imply that for every
$I=\{\zeta_1,\ldots ,\zeta_{k}\}$ size $k$ set of roots of $x^p=(-1)^{k+1}$ the object $T^{k(p-k)}_I$ is nonzero in the 
summand $\Fuk_\lambda(\Gr(k,p))$ of the Fukaya category, where $\lambda = pS_{\ydiagram{1}}(I) = p(\zeta_1 +\ldots +\zeta_{k})$.
\par
By Proposition \ref{Prop2} and Lemma \ref{DualSchur} $[\check{M}_{I^\vee}]\in T_{\cC^R}$ if and only if
$$S_{i\times j}(I) \neq 0 \quad \forall i\times j \quad \text{ rectangular diagram in } k\times (p-k) \text{ grid}\quad .$$
If by contradiction $S_{i\times j}(I)=0$ for some $I$ and diagram $i\times j$, then by definition
of Schur polynomial this means that
$$S_{i\times j}(I)=\sum_{T_{i\times j}}\zeta_1^{t_1}\cdots \zeta_{k}^{t_{k}}=0 \quad ,$$
where the sum runs over all semistandard Young tableaux $T_{i\times j}$ of the diagram $i\times j$,
as defined in Section \ref{Section1}.
Let us assume first that $k$ is odd, so that roots of $x^p=(-1)^{k+1}=1$ are $p$-th roots of unity.
Then each term in $S_{i\times j}(I)$ above is itself a $p$-th root of unity, and we have a vanishing
sum of a number of $p$-th roots of unity equal to $S_{i\times j}(1,\ldots ,1)$, i.e. the number
of semistandard Young tableaux of $i\times j$. A result of Lam-Leung on vanishing sums
of roots of unity \cite[Theorem 5.2]{LL} implies that $S_{i\times j}(1,\ldots ,1)$ must be a multiple
of $p$. On the other hand by Stanley's \textit{hook-content formula} (\cite[Theorem 15.3]{S} and Remark
\ref{RemarkThm3} below) we have
$$S_{i\times j}(1,\ldots ,1)=\prod_{u\in i\times j}\frac{k+c(u)}{h(u)}\quad .$$
The product ranges over boxes $u$ of the diagram $i\times j$. The content of
a box at entry $(s,t)$ of the grid is the number $c(u)=t-s$, and its hook number $h(u)$ is
the number of boxes of the diagram $i\times j$ below and to the right of $u$ (with $u$ itself counted
once).
\par
Therefore we must have
$$\prod_{u\in i\times j}\frac{k+c(u)}{h(u)}\equiv 0 \quad (\operatorname{mod}\, p)\quad .$$
By assumption $p$ is prime and $1\leq h(u)\leq p-1$, therefore there exists $u\in i\times j$ such that
$$c(u)\equiv -k \quad (\operatorname{mod}\, p)\quad .$$
Being $u\in i\times j$, it has to be at entry $(s,t)$ of the grid with $1\leq s\leq i$ and
$1\leq t \leq j$. Also notice that being the rectangle $i\times j$ in a $k\times (p-k)$ grid we have
$1\leq i\leq k$ and $1\leq j\leq p-k$. We conclude that
$$-k < 1-k\leq 1-i\leq c(u)=t-s \leq j-1\leq p-k-1 < p-k$$
in contradiction with $c(u)\equiv -k$ modulo $p$. This concludes the
proof in the case of $k$ odd. When $k$ is even $I$ consists of roots of $x^p=(-1)^{k+1}=-1$
and the argument above doesn't apply; on the other hand by Lemma \ref{DualSchur}
$$S_{i\times j}(I) \neq 0 \iff S_{j\times i}(e^{\pi i p}I^c) \neq 0$$
where $I^c$ denotes the roots of $x^p=(-1)^{k+1}=-1$ that are not in $I$, giving $p-k$ distinct roots
of $x^p=(-1)^{p-k+1}=1$ once rescaled by $e^{\pi i p}$. This reduces the problem to the previous
case and thus proves that the collection of $T^{k(p-k)}_I$ gives generators for all summands
of the Fukaya category.
\par
To prove homological mirror symmetry we argue as follows. Denoted $d=k(p-k)$, the assumption of $n=p$
prime guarantees that for every critical value $\lambda\in\mathbb{C}$ there is exactly one critical point
$[\check{M}_{I_\lambda^\vee}]\in\vGr(k,p)$ of $W$ with critical
value $\lambda$, and $[\check{M}_{I_\lambda^\vee}]\in T_{\cC^R}$ by the argument given earlier. Therefore
$[\check{M}_{I_\lambda^\vee}]=\theta_R(\hol_{I_\lambda})$ for a unique $\hol_{I_\lambda}\in(\mathbb{C}^\times)^d$ critical point
of $W_{T^d}$ with critical value $\lambda$. We denote $m_\lambda\subset\mathbb{C}[x_1^\pm,\ldots ,x_d^\pm]$
the maximal ideal corresponding to the point $\hol_{I_\lambda}$, and $T^d_{I_\lambda}$ the
generator of $\D\Fuk_\lambda(\Gr(k,p))$.
\par
The critical point $\hol_{I_\lambda}$ is nondegenerate. This condition holds because
a degenerate critical point of $W_{T^d}$ would correspond to a nonreduced point in the critical locus scheme
$Z\subset \vGr(k,p)$ of $W$ (see for example \cite[Lemma 3.5]{OT}), but closed mirror symmetry for Grassmannians (Theorem \ref{ClosedMS})
says that
$$Z=\operatorname{Spec}(\operatorname{Jac}(W))\cong\operatorname{Spec}(\QH(\Gr(k,p))) \quad ,$$
and this scheme is reduced because $\QH(\Gr(k,p))$ is semisimple, being
$$\QH(\Gr(k,p))=\bigoplus_{\lambda}\QH_\lambda(\Gr(k,p))$$
an algebra decomposition with one-dimensional summands.
From Theorem \ref{Fact4} of the Setup section we conclude that
$\D\Fuk_\lambda(\Gr(k,p))\simeq \D(\Cl_d)$, where $\Cl_d$ denotes the Clifford algebra of
the quadratic form of rank $d$ on $\mathbb{C}^d$.
Now combining the \textit{locality property} of the derived category of singularities established by
Orlov \cite[Proposition 1.14]{Or} and the fact proved above that all the critical points of $W$ are
in the rectangular chart $T_{\cC^R}\subset\vGr(k,p)$, we have for any $\lambda\in\mathbb{C}$
$$\D\mathcal{S}(W^{-1}(\lambda))\simeq \D\mathcal{S}(W^{-1}(\lambda)\cap T_{\cC^R})$$
where the intersection on the right is an affine scheme given by
$$W^{-1}(\lambda)\cap T_{\cC^R}=\operatorname{Spec}(\mathbb{C}[\check{p}_d^\pm: d \text{ rectangular}]/(W\restriction_{T_{\cC^R}}-\lambda))\quad .$$
Matrix factorizations are just another model for the category of singularity in the affine case, so that to conclude the proof of
homological mirror symmetry it suffices to establish an equivalence
$$\D(\Cl_d)\simeq \D\mathcal{M}(\mathbb{C}[\check{p}_d : d \text{ rectangular}],W\restriction_{T_{\cC^R}}-\lambda)\quad .$$
Dyckerhoff (Theorem 4.11 \cite{Dy}) shows that the localization ring morphism
$\mathbb{C}[x_1^\pm,\ldots ,x_d^\pm]\to \mathbb{C}[x_1^\pm,\ldots ,x_d^\pm]_{m_\lambda}$
induces an equivalence
$$\D\mathcal{M}(\mathbb{C}[x_1^\pm,\ldots ,x_d^\pm],W_{T^d}-\lambda)\simeq \D\mathcal{M}(\mathbb{C}[x_1^\pm,\ldots ,x_d^\pm]_{m_\lambda},W_{T^d}-\lambda)$$
and explicitly describes a generator \cite[Corollary 2.7]{Dy} of $\D\mathcal{M}(\mathbb{C}[x_1^\pm,\ldots ,x_d^\pm]_{m_\lambda},W_{T^d}-\lambda)$
whose endomorphism algebra is again the Clifford algebra $\Cl_d$ above, so that
$$\D\mathcal{M}(\mathbb{C}[x_1^\pm,\ldots ,x_d^\pm]_{m_\lambda},W_{T^d}-\lambda)\simeq\D(\Cl_d)$$
and this concludes the proof (see also Sheridan \cite[Section 6.1]{Sh} for a discussion of intrinsic
formality of Clifford algebras over $\mathbb{C}$).
\end{proof}

\begin{remark}\label{RemarkThm3}
The hook-content formula used in Theorem \ref{Thm3} can be obtained from \cite[Theorem 15.3]{S}
by setting $x=1$ in the generating function $H_m(\lambda)(x)$ of column-strict plane partitions
of shape $\lambda$ and largest part $\leq m$. In our application, $\lambda$ is the rectangular Young
diagram $i\times j$ in the $k\times (p-k)$ grid and $m=k$. In the terminology and notation
of Stanley's cited work, what we call semi-standard Young tableau on the Young diagram $\lambda$
in this article corresponds to a column-strict plane partition with shape $\lambda$ and
parts (or labels) in $S=\{1,\ldots ,k\}$. The correspondence is given by converting each label $l$ in
the Young tableau in the label $k+1-l$ of the corresponding plane partition,
so that the result has strictly decreasing columns and weakly decreasing rows. Also observe
that Stanley's definiton of Schur function $e_\lambda(x_1,\ldots ,x_k)$ in terms of plane
partitions with shape $\lambda$ agrees with our $S_\lambda(x_1,\ldots ,x_k)$ defined in terms of
tableaux on $\lambda$ because
$$S_\lambda(x_1,\ldots ,x_k)=e_\lambda(x_k,\ldots ,x_1)=e_\lambda(x_1,\ldots ,x_k)\quad ,$$
where the last equality holds because Schur polynomials are symmetric.
\end{remark}

\bibliographystyle{abbrv}
\bibliography{biblio}

\end{document}